\documentclass[leqno]{amsart}

\usepackage[foot]{amsaddr}

\usepackage{color}
\usepackage{geometry}
\usepackage{subfiles}
\usepackage{graphicx}
\usepackage[section]{placeins}
\usepackage{hyperref}
\usepackage{stmaryrd}
\usepackage{amsmath,amssymb,amsfonts,amsthm,textcomp}
\usepackage{mathtools}
\usepackage{tensor}
\usepackage{multicol}
\usepackage{subfig}
\usepackage{csquotes}
\usepackage{stackrel}
\usepackage{tikz}
\graphicspath{
              {Images/}
             }

\newfont{\tenbfsl}{cmbxti9 scaled 1200}
\newfont{\tenbbb}{msbm10}
\newfont{\svnbbb}{msbm8}

\DeclareMathOperator*{\argmin}{arg\,min}
\DeclareMathOperator*{\argmax}{arg\,max}

\newcommand{\bs}[1]{\boldsymbol{#1}}
\newcommand{\cl}[1]{\mathcal{#1}}

\newcommand{\bb}[1]{\mathbb{#1}}

\newcommand{\fr}[2]{\textstyle{\frac{{#1}}{{#2}}}}

\newcommand{\di}[1]{\,\mathrm{d}{#1}}

\newcommand{\trans}{\scriptscriptstyle\mskip-1mu\top\mskip-2mu}

\newcommand{\pprec}{\prec\mathrel{\mkern-7mu}\prec}

\newtheorem{rmk}{Remark}

\begin{document}

\title[Bayesian optimal experimental design with nuisance uncertainty]{Small-noise approximation for Bayesian optimal experimental design with nuisance uncertainty}
\author{Arved Bartuska$^1$, Luis Espath$^2$ \& Ra\'{u}l Tempone$^{1,3,4}$}
\address{$^1$Department of Mathematics, RWTH Aachen University, Geb\"{a}ude-1953 1.OG, Pontdriesch 14-16, 161, 52062 Aachen, Germany.}
\address{$^2$School of Mathematical Sciences, University of Nottingham, Nottingham, NG7 2RD, United Kingdom.}
\address{$^3$King Abdullah University of Science \& Technology (KAUST), Computer, Electrical and Mathematical Sciences \& Engineering Division (CEMSE), Thuwal 23955-6900, Saudi Arabia.}
\address{$^4$Alexander von Humboldt Professor in Mathematics for Uncertainty Quantification, RWTH Aachen University, Germany.}
\email{bartuska@uq.rwth-aachen.de}

\date{\today}

\begin{abstract}
\noindent
Calculating the expected information gain in optimal Bayesian experimental design
typically relies on nested Monte Carlo sampling. When the model also contains nuisance parameters, which are parameters that contribute to the overall uncertainty of the system but are of no interest in the Bayesian design framework, this introduces
a second inner loop. We propose and derive a small-noise approximation for
this additional inner loop.
The computational cost of our method can be further reduced by applying a
Laplace approximation to the remaining inner loop. Thus, we present two methods, the small-noise double-loop Monte Carlo and small-noise Monte Carlo Laplace methods.
Moreover, we demonstrate that the total complexity of these two approaches remains comparable to the case without nuisance uncertainty. To assess the efficiency of these methods, we present three examples, and the last example includes the partial differential equation for the electrical impedance tomography experiment for composite laminate materials.
\\
\textbf{AMS subject classifications:}
$\cdot$
65K05 
$\cdot$
65C05 
$\cdot$

\end{abstract}

\maketitle

\tableofcontents                        


\section{Introduction}
When conducting scientific experiments, information from data is used to learn about the parameters of an underlying model \cite{Cha95}. The relationship between data and parameters depends on the experiment design, which can be regarded as a hyperparameter. Data acquisition is often the most expensive part of an experiment; thus, this design parameter can be optimized beforehand to reduce the number of data points required to achieve a certain error tolerance in the parameters. The amount of information obtained from an experiment can be expressed as the expected Kullback--Leibler divergence \cite{Kul51, Kul59}, also known as the expected information gain (EIG) \cite{Lin56}.

Computing this quantity requires evaluating two nested integrals, which are typically not given in closed form. As the straightforward approximation of these integrals using the Monte Carlo method \cite{Rya03} is often prohibitively expensive, several more advanced estimators based on the Laplace approximation \cite{Sti86, Tie86, Tie89, Kas90} have been proposed \cite{Bec18, Bec20, Car20, Fen19}. These estimators either approximate the inner integral using a second-order Taylor expansion around the mode of the posterior probability density function (pdf) directly or use this approximation for importance sampling to reduce the number of samples required for the inner Monte Carlo estimation. The first approach reduces the order of the average work required but introduces a bias that depends on the number of repetitive experiments, while the second approach has a controllable bias but only reduces constant terms in the asymptotic work requirement. For more details on the Laplace approximation, see \cite{Sch20, Wac17}, where a bound on the approximation error in terms of the Hellinger distance is derived. The approximation error of Laplace-based importance sampling is detailed in \cite{Sch20}. These are asymptotic results for the small noise or large data limit. We only consider a fixed noise level and fixed amount of data. Recently, \cite{Hel22} derived pre-asymptotic error estimates for the Laplace approximation and investigate the effects of the nonlinearity of the forward problem on the non-Gaussianity of the likelihood. Furthermore, \cite{Spo22} provide non-asymptotic guarantees for the Laplace approximation in total variation distance and also consider the case where the maximum aposteriori estimate and Hessian of the log-likelihood are not available directly, but have to be approximated. The concentration of the posterior pdf for time series data and the resulting increased accuracy of the Laplace approximation is discussed in \cite{Lon15b}.
If the prior pdf has compact support, e.g., a uniform distribution, the posterior distribution exists on the same bounded domain. The Gaussian distribution resulting from the Laplace approximation, however, is unbounded and thus introduces an additional error. In \cite{Bis16}, a truncated version of the Laplace approximation is introduced alongside a rejection sampling algorithm to efficiently estimate the expected information gain under this setting.

If the model also contains unknown parameters that are not of immediate interest, we refer to them as nuisance parameters, and the usual approach is to marginalize them \cite{Fen19, Ber79, Pol88, Lie13}, leading to a second inner integral and increasing the computational burden. Such a scenario may arise if some of the model's parameters are not well-calibrated. We will show this for a practical example where the electric conductivity of a material is only given by a concentrated log-normal distribution rather than an exact value. Assuming that the effect of the nuisance parameters is small, we propose using a small-noise approximation combined with the Laplace approximation. The small-noise approximation uses a first-order Taylor approximation of the residual term and a first-order approximation of the resulting likelihood. This method works if the nuisance covariance is sufficiently small compared to the covariance of the error to ensure positive definiteness of the resulting covariance matrix. A slightly different approach was used in \cite{Lon21}, where the error covariance was considered an unknown parameter, also requiring marginalization. Another approach is followed in \cite{Ken01}, where the authors propose to estimate hyperparameters beforehand and then consider them to be fixed to avoid high-dimensional quadrature for the marginalization.

Our approach differs from that in \cite{Fen19} in the estimators proposed to approximate the EIG functional. Both \cite{Fen19}  and our work start out with the same nuisance parameter model, requiring double-loop Monte Carlo approximation with two inner loops. Our contribution lies in proposing a novel small-noise approximation, thus avoiding this additional inner loop at the cost of a bound on the maximum allowed nuisance uncertainty. This bound is consistent with our differing understanding of the interpretation of nuisance parameters. A comparison of the computational work is provided by Remark 1, where we show the optimal number of required samples for two inner loops. We use the small-noise approximation to handle the nuisance uncertainty in a novel manner, which allows us to construct two estimators: one based on the Laplace approximation for the inner integral and one using importance sampling based on the Laplace approximation. The first estimator consists of just one loop but introduces a fixed bias, whereas the second estimator consists of one outer and one inner loop and introduces a controllable bias. The work \cite{Fen19} constructs a multiple-importance sampling scheme that reuses previous model evaluations and consists of one outer and two inner loops. There is also a difference in how we justify the need to deal with nuisance parameters in this model. We aim to extend the data model by including nuisance parameters. In contrast, \cite{Fen19} aimed to reduce the parameter space dimension by regarding some existing parameters as nuisance parameters. Following the distinction between aleatoric and epistemic uncertainty in \cite{Ale20}, this means that \cite{Fen19} viewed the nuisance parameters as an epistemic uncertainty and assumed that something could in principle be learned about them. Our approach is agnostic as to the type of uncertainty. The nuisance parameters could also be viewed as an aleatoric uncertainty inaccessible to Bayesian updates, as nothing can be learned about them independent of the chosen design, even for an unlimited amount of data. Similar to \cite{Fen19, Lev14}, we consider the distinction between parameters of interest and nuisance parameters to be a modeling choice. As \cite{Ber79} points out, the marginalization of nuisance parameters can also be formulated in terms of a Quantity of Interest (QoI). The application of the Laplace approximation to optimal design problems w.r.t. a QoI was discussed in \cite{Lon13, Lon15}.

In this paper, we introduce the Bayesian setting in Section \ref{sec:bayesian.formulation} and develop the small-noise approximation in Section \ref{sec:small.noise.approximation}. We use the resulting updated likelihood function in Section \ref{sec:double.loop} and formulate the double-loop Monte Carlo (DLMC) estimator along with the optimal number of outer and inner samples and the error-splitting term. We formulate the Monte Carlo Laplace (MCLA) estimator based on the small-noise approximation in Section \ref{sec:laplace} and the DLMC estimator with importance sampling (DLMCIS) in Section \ref{sec:importance.sampling} and derive the optimal setting for both estimators. In Section \ref{sec:numerical.results}, we present three numerical examples, one to showcase the degradation of the small-noise approximation for a large nuisance error, one linear Gaussian example as presented in \cite{Fen19} to determine the effect of nuisance parameters on the optimal design, and one based on the partial differential equation (PDE) for electrical impedance tomography (EIT) to establish the practical use of the proposed method.

\section{Bayesian formulation}\label{sec:bayesian.formulation}

Our goal is to determine a design that lets us optimally estimate a vector-valued parameter of interest $\bs{\theta}$ from measured data $\bs{y}$ given a mathematical model $\bs{g}$ of an experiment depending on $\bs{\theta}$, $\bs{y}$, and a design $\bs{\xi}$. In the EIT example \cite{Som92} discussed in Section \ref{sec:EIT}, we consider the fiber orientation of laminate materials as our parameter of interest and the locations of measurement electrodes as the design.

We model the data with additive noise $\bs{\epsilon}$, that is,
\begin{equation}\label{eq:data.model}
\bs{y}_i(\bs{\xi})=\bs{g}(\bs{\xi},\text{$\bs{\theta}_t$},\text{$\bs{\phi}_t$})+\text{$\bs{\epsilon}_i$}.
\end{equation}
If the model also contains additional parameters $\bs{\phi}$ that are of no interest in the Bayesian inference formulation, we refer to them as nuisance parameters and marginalize them. This method can be computationally prohibitive; therefore, other methods are addressed to deal with nuisance parameters. In the EIT example, electric conductivities of the materials are considered nuisance parameters. The model contains the following quantities
\begin{itemize}
  \item $\bs{y}_i\in\bb{R}^q$ is the vector of measurements, $1\leq i\leq N_e$;
  \item $\bs{Y}=(\bs{y}_1,\ldots,\bs{y}_{N_e}) \in \bb{R}^{q\times N_e}$ are the measurements for each experiment;
  \item $\bs{g}\in\bb{R}^q$ denotes the model predictions;
  \item $\bs{\epsilon}_i\in\bb{R}^q$ is the error vector assumed to be Gaussian $\bs{\epsilon}\sim\cl{N}(\bs{0},\bs{\Sigma}_{\bs{\epsilon}})$;
  \item $\bs{\theta}_t\in\bb{R}^{d_{\bs{\theta}}}$ represents the true parameter vector;
  \item $\bs{\phi}_t\in\bb{R}^{d_{\bs{\phi}}}$ denotes the true vector of the nuisance parameters;
	\item $\bs{\xi}\in\bb{R}^s$ is the vector of design parameters;
	\item $q$ indicates the number of measurements; and
	\item $N_e$ is the number of observations under consistent experimental setup.
\end{itemize}
In this context, the Bayes formula reads
\begin{align}\label{eq:posterior}
\pi(\bs{\theta}|\bs{Y},\bs\xi) &=
 \dfrac{\pi(\bs{\theta}) p(\bs{Y}|\bs{\theta},\bs\xi)}{p(\bs{Y}|\bs\xi)}\nonumber\\
 &=\dfrac{\pi(\bs{\theta})\int_{\bs{\Phi}} p(\bs{Y}|\bs{\theta},\bs{\phi},\bs\xi)\pi(\bs{\phi})\di{}\bs{\phi}}{\int_{\bs{\Theta}}\int_{\bs{\Phi}}p(\bs{Y}|\bs{\theta},\bs{\phi},\bs\xi)\pi(\bs{\phi})\pi(\bs{\theta})\di{}\bs{\phi}\di{}\bs{\theta}},
\end{align}
where the likelihood $p(\bs{Y}|\bs{\theta},\bs\xi)$ and evidence $p(\bs{Y}|\bs\xi)$ terms are obtained by marginalization. We distinguish between $p(\cdot)$ for pdfs of the data and $\pi(\cdot)$ for pdfs of model parameters, namely, the prior distribution $\pi(\bs{\theta})$ and the posterior distribution $\pi(\bs{\theta}|\bs{Y},\bs\xi)$. The former signifies our knowledge before conducting any experiments, while the latter signifies our knowledge once data is obtained. The term $\pi(\bs{\phi})$ is a regular probability density and not a prior distribution that is updated.

Bayesian optimal experimental design aims to determine the experimental setup, as described by the design parameter $\bs{\xi}$, for the Bayesian inference of $\bs{\theta}_t$. We consider batch design similar to \cite{Lon21}, where different experimental setups can be expressed by the parameter $\bs{\xi}$, which may be a boolean-, integer-, real-, or complex-valued scalar or vector variable. This parameter depends on the particular physical experiment. To carry out the design optimization, we need to establish a relation between an experimental setup and the amount of information that it provides. The information gain based on the Shannon entropy \cite{Sha48} for an experimental design $\bs{\xi}$ is commonly measured by the Kullback--Leibler divergence \cite{Kul51, Kul59}. The Kullback--Leibler divergence is given by
\begin{equation}\label{Dkl}
  D_{\mathrm{KL}}=\int_{\bs{\Theta}}[\log(\pi(\bs{\theta}|\bs{Y}))-\log(\pi(\bs{\theta}))]\pi(\bs{\theta}|\bs{Y})\di{}\bs{\theta},
\end{equation}
where the design parameter $\bs\xi$ is omitted for conciseness. The setup must be determined before any experiments are conducted; thus, the marginalization over the data $\bs{Y}$ must be considered, which is referred to as the EIG \cite{Lin56}:
\begin{equation}\label{eq:EIG.post}
  I=\int_{\cl{Y}}\int_{\bs{\Theta}}[\log(\pi(\bs{\theta}|\bs{Y}))-\log(\pi(\bs{\theta}))]\pi(\bs{\theta}|\bs{Y})\di{}\bs{\theta} p(\bs{Y})\di{}\bs{Y}.
\end{equation}
We replace the posterior pdf using the Bayes formula \eqref{eq:posterior} to arrive at:
\begin{equation}\label{eq:EIG.lik}
  I=\int_{\bs{\Theta}}\int_{\cl{Y}}[\log(p(\bs{Y}|\bs{\theta}))-\log(p(\bs{Y}))]p(\bs{Y}|\bs{\theta}) \di{}\bs{Y}\pi(\bs{\theta})\di{}\bs{\theta}.
\end{equation}
The goal is to learn about the parameter of interest $\bs{\theta}$ and not the nuisance parameter $\bs{\phi}$.

The likelihood function is given as follows:
\begin{equation}\label{eq:likelihood}
p(\bs{Y}|\bs{\theta},\bs{\phi})\coloneqq\det(2\pi\bs{\Sigma}_{\bs{\varepsilon}})^{-\frac{N_e}{2}}\exp\left(-\frac{1}{2}\sum_{i=1}^{N_e} \bs{r}_i(\bs{\theta},\bs{\phi},\bs{\epsilon}_i) \cdot \bs{\Sigma}_{\bs{\varepsilon}}^{-1}\bs{r}_i(\bs{\theta},\bs{\phi},\bs{\epsilon}_i)\right),
\end{equation}
where
\begin{equation}\label{eq:residual}
\bs{r}_i(\bs{\theta},\bs{\phi},\bs{\epsilon}_i)\coloneqq\bs{y}_i-\bs{g}(\bs{\theta},\bs{\phi})=\bs{g}(\bs{\theta}_t,\bs{\phi}_t)+\bs{\epsilon}_i-\bs{g}(\bs{\theta},\bs{\phi}).
\end{equation}

In view of \eqref{eq:posterior}, the EIG \eqref{eq:EIG.lik} reads
\begin{align}\label{eq:EIG.definition}
  I&=\int_{\bs{\Theta}}\int_{\bs{\Phi}}\int_{\cl{Y}}\Bigg[\log\left(\int_{\bs{\Phi}} p(\bs{Y}|\bs{\theta},\bs{\varphi})\pi(\bs{\varphi})\di{}\bs{\varphi}\right)\nonumber\\[4pt]
  &-\log\left(\int_{\bs{\Theta}}\int_{\bs{\Phi}}p(\bs{Y}|\bs{\vartheta},\bs{\varphi})\pi(\bs{\varphi})\pi(\bs{\vartheta})\di{}\bs{\varphi}\di{}\bs{\vartheta}\right)\Bigg]p(\bs{Y}|\bs{\theta},\bs{\phi})\di{}\bs{Y}\pi(\bs{\phi})\di{}\bs{\phi}\pi(\bs{\theta}) \di{}\bs{\theta},
\end{align}
where $\bs\varphi$ and $\bs\vartheta$ are dummy variables (to distinguish them from $\bs{\phi}$ and $\bs{\theta}$ from the outer integral). The expression \eqref{eq:EIG.definition} is typically not given in closed-form and estimators need to be chosen carefully. In a sample-based version of \eqref{eq:EIG.definition}, we arrive at the DLMC estimator with two inner loops
\begin{equation}\label{EIG}
 \frac{1}{N}\sum_{n=1}^N\left[\log\left(\frac{1}{M_1}\sum_{m=1}^{M_1}p(\bs{Y}^{(n)}|\bs{\theta}^{(n)},\bs{\varphi}^{(n,m)})\right)-\log\left(\frac{1}{M_2}\sum_{k=1}^{M_2}p(\bs{Y}^{(n)}|\bs{\vartheta}^{(n,k)},\bs{\varphi}^{(n,k)})\right)\right],
\end{equation}
where samples are drawn as follows. We draw $\bs{\theta}^{(n)}\stackrel{\mathrm{iid}}{\sim}\pi(\bs{\theta})$, $\bs{\phi}^{(n)}\stackrel{\mathrm{iid}}{\sim}\pi(\bs{\phi})$, $\bs{Y}^{(n)}\stackrel{\mathrm{iid}}{\sim}p(\bs{Y}|\bs{\theta}^{(n)},\bs{\phi}^{(n)})$, and $\bs{\varphi}^{(n,m)}\stackrel{\mathrm{iid}}{\sim}\pi(\bs{\varphi})$ for the first term and $\bs{\vartheta}^{(n,k)}\stackrel{\mathrm{iid}}{\sim}\pi(\bs{\theta})$ and $\bs{\varphi}^{(n,k)}\stackrel{\mathrm{iid}}{\sim}\pi(\bs{\phi})$ for the second term. To make clear the difference between the EIG formulation considered in this paper and others, e.g., \cite{Rya03, Bec18, Lon13}, we make the following remark.

\begin{rmk}[Differences arising from nuisance parameters]\label{rk:diff}
Marginalization is necessary for the denominator and numerator of the log-likelihood ratio, compared to the nuisance-free formulation \cite{Bec18}:
 \begin{equation}
   I=\int_\Theta\int_\cl{Y}\left[\log(p(\bs{Y}|\bs{\theta}))-\log\left(\int_\Theta p(\bs{Y}|\bs{\vartheta})\pi(\bs{\vartheta})\di{}\bs{\vartheta}\right)\right]p(\bs{Y}|\bs{\theta})\di{}\bs{Y}\pi(\bs{\theta})\di{}\bs{\theta}
 \end{equation}
with the corresponding estimator
\begin{equation}
  \frac{1}{N}\sum_{n=1}^N\left[\log\left(p(\bs{Y}^{(n)}|\bs{\theta}^{(n)})\right)-\log\left(\frac{1}{M}\sum_{m=1}^{M}p(\bs{Y}^{(n)}|\bs{\vartheta}^{(n,m)})\right)\right].
\end{equation}
Instead of one outer and one inner loop, we now have two inner loops. The cost associated with performing the first inner loop in \eqref{EIG} is proportional to the conditional variance of the likelihood $M_1\propto\bb{V}[p(\bs{Y}|\bs{\theta},\bs{\phi})|\bs{Y},\bs{\theta}]\eqqcolon C_1$ and the cost associated with the second inner loop is proportional to $M_2\propto\bb{V}[p(\bs{Y}|\bs{\theta},\bs{\phi})|\bs{Y}]\eqqcolon C_2$. From the optimal setting derivation in this case, it follows that (see Appendix \ref{ap:double.inner.loop})
\begin{equation}
M_1+M_2=\frac{C_1+2\sqrt{C_1C_2}+C_2}{\text{const.}}.
\end{equation}
For $C_2=0$, this gives the work of the double-loop estimator with only one inner loop, but for $C_1=C_2$, this increases the work by a factor of four. This also means that for $C_2=(3-2\sqrt{2})C_1\approx0.17C_1$, we have an increase of a factor of two. Next, we derive the small-noise approximation that lets us carry out the marginalization of the nuisance parameters without introducing an additional inner loop.
\end{rmk}

\section{Small-noise approximation}\label{sec:small.noise.approximation}

Any Gaussian random vector may be expressed as the rotation of a random vector with independent components. Thus, without loss of generality, let $\bs{\phi}=\bs{\Sigma}^{1/2}_{\bs{\phi}} \bs{z}$, where the components of $\bs{z}\sim \cl{N}(\bs{0},\bs{1})$ are independent and $\bs{1}$ is the identity matrix with the appropriate dimension. We further assume that $\bs{\Sigma}^{1/2}_{\bs{\phi}}\preceq\bs{1}$. Next, the Taylor expansion on \eqref{eq:residual} reads
\begin{equation}\label{eq:small.residual}
  \bs{r}_i(\bs{\theta},\bs{\phi},\bs{\epsilon}_i)=\bs{r}_i(\bs{\theta},\bs{0},\bs{\epsilon}_i)-\nabla_{\bs{\phi}}\bs{g}(\bs{\theta},\bs{0})\bs{\Sigma}^{1/2}_{\bs{\phi}} \bs{z}+\cl{O}_{\bb{P}}\left(\|\bs{\Sigma}^{1/2}_{\bs{\phi}} \bs{z}\|^2\right),
\end{equation}
where $\bs{r}_i(\bs{\theta},\bs{0},\bs{\epsilon}_i)$ is the residual evaluated without nuisance parameters, $\bs{r}_i(\bs{\theta},\bs{\phi}=\bs{0},\bs{\epsilon}_i)$, for the $i^{\text{th}}$ experiment, where $1\leq i\leq N_e$. We also abuse notation and write $\cl{O}_{\bb{P}}\left(\|\bs{\Sigma}^{1/2}_{\bs{\phi}} \bs{z}\|^2\right)$ in a component-wise sense.\footnote{The notation $X_M=\cl{O}_{\bb{P}}(a_M)$ for a sequence of random variables $X_M$ and constants $a_M$ is as follows. For any $\epsilon>0$, there exists a finite $K(\epsilon)>0$ and finite $M_0>0$ such that $\bb{P}(|X_M|>K(\epsilon)|a_M|)<\epsilon$ holds for all $M\geq M_0$.}
We obtain the following expression for the likelihood \eqref{eq:likelihood}:
\begin{align}\label{eq:marginalization.likelihood}
  p(\bs{Y}|\bs{\theta})&=\det(2\pi\bs{\Sigma}_{\bs\epsilon})^{-\frac{N_e}{2}}\int_{\bs{\Phi}}\prod_{i=1}^{N_e}\exp\left(-\frac{1}{2} \bs{r}_i(\bs{\theta},\bs{\phi},\bs{\epsilon}_i)\cdot\bs{\Sigma}_{\bs\epsilon}^{-1}\bs{r}_i(\bs{\theta},\bs{\phi},\bs{\epsilon}_i) \right)\pi(\bs{\phi})\di{}\bs{\phi}\nonumber\\
  &=\det(2\pi\bs{\Sigma}_{\bs\epsilon})^{-\frac{N_e}{2}}\prod_{i=1}^{N_e}\bb{E}\left[\exp\left(-\frac{1}{2} \bs{r}_i(\bs{\theta},\bs{\phi},\bs{\epsilon}_i)\cdot\bs{\Sigma}_{\epsilon}^{-1}\bs{r}_i(\bs{\theta},\bs{\phi},\bs{\epsilon}_i)\right)\Big|\bs{\theta},\bs{\epsilon}_i\right].
\end{align}
The expectation in \eqref{eq:marginalization.likelihood} can be simplified using \eqref{eq:small.residual} as follows:
\begin{multline}\label{eq:small.exp}
  \bb{E}\left[\exp\left(-\frac{1}{2}\bs{r}_i(\bs{\theta},\bs{\phi},\bs{\epsilon}_i)\cdot\bs{\Sigma}_{\bs\epsilon}^{-1}\bs{r}_i(\bs{\theta},\bs{\phi},\bs{\epsilon}_i)\right)\Big|\bs{\theta},\bs{\epsilon}_i\right]\\
  =\exp\left(-\frac{1}{2}\bs{r}_i(\bs{\theta},\bs{0},\bs{\epsilon}_i)\cdot\bs{\Sigma}_{\bs\epsilon}^{-1}\bs{r}_i(\bs{\theta},\bs{0},\bs{\epsilon}_i)\right)\bb{E}\left[\exp\left(\bs{r}_i(\bs{\theta},\bs{0},\bs{\epsilon}_i)\cdot\bs{\Sigma}_{\bs\epsilon}^{-1}
  (\nabla_{\bs{\phi}}\bs{g}(\bs{\theta},\bs{0})\bs{\Sigma}^{1/2}_{\bs{\phi}}\bs{z})+\cl{O}_{\bb{P}}\left(\|\bs{\Sigma}^{1/2}_{\bs{\phi}} \bs{z}\|^2\right)\right)\Big|\bs{\theta},\bs{\epsilon}_i\right].
\end{multline}
Notice that the higher-order term is different from the one in \eqref{eq:small.residual}. Letting $\bs{w}\coloneqq\bs{r}_i(\bs{\theta},\bs{0},\bs{\epsilon}_i)\cdot\bs{\Sigma}_{\bs\epsilon}^{-1}
\nabla_{\bs{\phi}}\bs{g}(\bs{\theta},\bs{0})\bs{\Sigma}^{1/2}_{\bs{\phi}}$, we have that
\begin{align}\label{eq:small.exp.simplified}
  &\bb{E}\left[\exp\left(-\frac{1}{2}\bs{r}_i(\bs{\theta},\bs{\phi},\bs{\epsilon}_i)\cdot\bs{\Sigma}_{\bs\epsilon}^{-1}\bs{r}_i(\bs{\theta},\bs{\phi},\bs{\epsilon}_i)\right)\Big|\bs{\theta},\bs{\epsilon}_i\right]\nonumber\\
  &=\exp\left(-\frac{1}{2}\bs{r}_i(\bs{\theta},\bs{0},\bs{\epsilon}_i)\cdot\bs{\Sigma}_{\bs\epsilon}^{-1}\bs{r}_i(\bs{\theta},\bs{0},\bs{\epsilon}_i)\right)\bb{E}\left[\exp\left(\bs{w}\cdot\bs{z}+\cl{O}_{\bb{P}}\left(\|\bs{\Sigma}^{1/2}_{\bs{\phi}} \bs{z}\|^2\right)\right)\Big|\bs{\theta},\bs{\epsilon}_i\right],\nonumber\\
  &=\exp\left(-\frac{1}{2}\bs{r}_i(\bs{\theta},\bs{0},\bs{\epsilon}_i)\cdot\bs{\Sigma}_{\bs\epsilon}^{-1}\bs{r}_i(\bs{\theta},\bs{0},\bs{\epsilon}_i)\right)\left(\exp\left(\fr{1}{2}\bs{w}\cdot\bs{w}\right)+\bb{E}\left[\cl{O}_{\bb{P}}\left(\|\bs{\Sigma}^{1/2}_{\bs{\phi}} \bs{z}\|^2\right)\Big|\bs{\theta},\bs{\epsilon}_i\right]\right),
\end{align}
where we use that $\bb{E}[e^{tZ}]=e^{\frac{t^2}{2}}$ for $Z\sim\cl{N}(0,1)$.
Using tensor notation, we write the first term in \eqref{eq:small.exp.simplified} as
\begin{align}
  \bs{r}_i(\bs{\theta},\bs{0},\bs{\epsilon}_i)\cdot\bs{\Sigma}_{\bs\epsilon}^{-1}\bs{r}_i(\bs{\theta},\bs{0},\bs{\epsilon}_i)=\bs{r}_i(\bs{\theta},\bs{0},\bs{\epsilon}_i)\otimes\bs{\Sigma}_{\bs\epsilon}^{-1}\bs{r}_i(\bs{\theta},\bs{0},\bs{\epsilon}_i)\colon\bs{1},
\end{align}
and the second term as
\begin{equation}
  \bs{r}_i(\bs{\theta},\bs{0},\bs{\epsilon}_i)\cdot\bs{\Sigma}_{\bs\epsilon}^{-1}\bs{M}\bs{M}^{\trans}\bs{\Sigma}_{\bs\epsilon}^{-1}\bs{r}_i(\bs{\theta},\bs{0},\bs{\epsilon}_i)=\bs{r}_i(\bs{\theta},\bs{0},\bs{\epsilon}_i)\otimes\bs{\Sigma}_{\bs\epsilon}^{-1}\bs{r}_i(\bs{\theta},\bs{0},\bs{\epsilon}_i)\colon \bs{\Sigma}_{\bs\epsilon}^{-1}\bs{M} \bs{M}^{\trans},
\end{equation}
where $\bs{M}\coloneqq\nabla_{\bs{\phi}}\bs{g}(\bs{\theta},\bs{0})\bs{\Sigma}^{1/2}_{\bs{\phi}}$.
Thus, the full likelihood is written as follows:
\begin{equation}\label{eq:small.likelihood}
  p(\bs{Y}|\bs{\theta})\propto\exp\Bigg(-\frac{1}{2}\sum_{i=1}^{N_e}\bs{r}_i(\bs{\theta},\bs{0},\bs{\epsilon}_i)\otimes\bs{\Sigma}_{\bs\epsilon}^{-1}\bs{r}_i(\bs{\theta},\bs{0},\bs{\epsilon}_i)\colon\left(\bs{1}-\bs{\Sigma}_{\bs\epsilon}^{-1}\bs{M}\bs{M}^{\trans}\right)  +\cl{O}_{\bb{P}}\left(\|\bs{\Sigma}^{1/2}_{\bs{\phi}} \bs{z}\|^2\right)\Bigg),
\end{equation}
or in a more familiar form
\begin{multline}\label{eq:small.likelihood.fam}
  p(\bs{Y}|\bs{\theta})=\det\left(2\pi\Big[\bs{\Sigma}_{\bs\epsilon}^{-1}-\bs{\Sigma}_{\bs\epsilon}^{-1}\bs{M}\bs{M}^{\trans}\bs{\Sigma}_{\bs\epsilon}^{-1}\Big]^{-1}\right)^{-\frac{N_e}{2}}\\\exp\Bigg(-\frac{1}{2}\sum_{i=1}^{N_e}\bs{r}_i(\bs{\theta},\bs{0},\bs{\epsilon}_i)\cdot\Big[\bs{\Sigma}_{\bs\epsilon}^{-1}-\bs{\Sigma}_{\bs\epsilon}^{-1}\bs{M}\bs{M}^{\trans}\bs{\Sigma}_{\bs\epsilon}^{-1}\Big]\bs{r}_i(\bs{\theta},\bs{0},\bs{\epsilon}_i)
  +\cl{O}_{\bb{P}}\left(\|\bs{\Sigma}^{1/2}_{\bs{\phi}} \bs{z}\|^2\right)\Bigg),
\end{multline}
where we obtain the correction term
\begin{equation}\label{eq:correction.term}
\bs{C}\coloneqq\bs{1}-\bs{\Sigma}_{\bs\epsilon}^{-1/2}\bs{M}\bs{M}^{\trans}\bs{\Sigma}_{\bs\epsilon}^{-1/2},
\end{equation}
depending on $\nabla_{\bs{\phi}}\bs{g}$ and (small) $\bs{\Sigma}_{\bs{\phi}}$. As $\bs{\Sigma}_{\bs{\phi}}$ increases, this term becomes smaller, which has the same effect as $\bs{\Sigma}_\epsilon$ increasing. We will now provide more details on the applicability of the small-noise approximation in the following remarks.

\begin{rmk}[Conditions for correction term $\bs{C}$]\label{rmk.eig}
[$\bs{\Sigma}_{\bs\epsilon}^{-1/2}\bs{C}\bs{\Sigma}_{\bs\epsilon}^{-1/2}$ as a metric tensor.] To ensure that \eqref{eq:small.likelihood.fam} can be understood as a distance function, $\bs{\Sigma}_{\bs\epsilon}^{-1/2}\bs{C}\bs{\Sigma}_{\bs\epsilon}^{-1/2}$ must be a metric tensor. That is, it must be positive definite, which is only the case if $\bs{\Sigma}_{\bs{\phi}}^{1/2}\pprec\bs{\Sigma}_{\bs{\epsilon}}^{1/2}$.\footnote{The symbol in $\bs{A}\pprec\bs{B}$ means that the matrix $\bs{A}$ has much smaller eigenvalues than $\bs{B}$.}
\end{rmk}

\begin{rmk}[Special case for $\bs{\Sigma}=\sigma^2 \bs{1}$]\label{rmk:spc}
For the special case in which $\bs{\Sigma}_{\bs{\epsilon}}=\sigma_{\bs{\epsilon}}^2 \bs{1}$ and $\bs{\Sigma}_{\bs{\phi}}=\sigma_{\bs{\phi}}^2 \bs{1}$, we make the previous result more specific. We rewrite the matrix $\bs{C}$ as
\begin{align}\label{eq:eig.Jac}
\bs{C}&=\bs{1}-\bs{\Sigma}_{\bs{\epsilon}}^{-1/2}\bs{M}\bs{M}^{\trans}\bs{\Sigma}_{\bs{\epsilon}}^{-1/2}=\bs{1}-\frac{\sigma_{\bs{\phi}}^2}{\sigma_{\bs{\epsilon}}^2}\nabla_{\bs{\phi}}\bs{g}(\bs{\theta},\bs{0})\nabla_{\bs{\phi}}\bs{g}(\bs{\theta},\bs{0})^{\trans},
\end{align}
which is positive definite if all its eigenvalues $\mu_i$, $1\leq i\leq q$, are positive. These eigenvalues can be written as
\begin{equation}
\mu_i=1-\frac{\sigma_{\bs{\phi}}^2}{\sigma_{\bs{\epsilon}}^2}\lambda_i, \quad 1\leq i\leq q,
\end{equation}
where $\lambda_i$ are the eigenvalues of $\nabla_{\bs{\phi}}\bs{g}(\bs{\theta},\bs{0})\nabla_{\bs{\phi}}\bs{g}(\bs{\theta},\bs{0})^{\trans}$. It can easily be seen that this matrix has rank 1 and that the only nonzero eigenvalue is given by $\lVert \nabla_{\bs{\phi}}\bs{g}(\bs{\theta},\bs{0})\rVert_2^2$. This quantity is a random variable dependent on $\bs{\theta}$, so
\begin{equation}
  \sigma_{\bs{\phi}}^2\lVert \nabla_{\bs{\phi}}\bs{g}(\bs{\theta},\bs{0})\rVert_2^2<\sigma_{\bs{\epsilon}}^2
\end{equation}
 must hold almost surely to guarantee positive definiteness of the updated covariance matrix $\bs{\Sigma}_{\bs\epsilon}^{-1/2}\bs{C}\bs{\Sigma}_{\bs\epsilon}^{-1/2}$. So either the gradient of $\bs{g}$ or the covariance $\sigma_{\bs{\phi}}^2$ have to be small. Thus, either the nuisance parameters have little impact on the model or there is only a small uncertainty coming from the nuisance parameters.
\end{rmk}

\begin{rmk}[Small-noise approximation for linear models]
For linear functions on $\bs{\phi}$, the approximation \eqref{eq:small.residual} is exact, as $\nabla_{\bs{\phi}}\nabla_{\bs{\phi}}\bs{g}=\bs{0}$, but the second approximation \eqref{eq:small.exp} is not exact, and one must still ensure that the updated covariance matrix $\bs{\Sigma}_{\bs\epsilon}^{-1/2}\bs{C}\bs{\Sigma}_{\bs\epsilon}^{-1/2}$ remains positive definite.
\end{rmk}
We can now employ the small-noise approximation in several numerical estimators for the EIG \eqref{eq:EIG.definition}.

\section{Double-loop Monte Carlo}\label{sec:double.loop}
\subsection{Double-loop Monte Carlo estimator}
We obtained an approximation of the likelihood \eqref{eq:likelihood} given via \eqref{eq:small.likelihood.fam} as
\begin{equation}\label{eq:likelihood.tilde}
  \tilde{p}(\bs{Y}|\bs{\theta})=\det\left(2\pi\bs{\Sigma}_{\bs\epsilon}^{1/2}\bs{C}^{-1}\bs{\Sigma}_{\bs\epsilon}^{1/2}\right)^{-\frac{N_e}{2}}\exp\left(-\frac{1}{2}\sum_{i=1}^{N_e}\bs{r}_i(\bs{\theta},\bs{0},\bs{\epsilon}_i)\cdot\bs{\Sigma}_{\bs\epsilon}^{-1/2}\bs{C}\bs{\Sigma}_{\bs\epsilon}^{-1/2}\bs{r}_i(\bs{\theta},\bs{0},\bs{\epsilon}_i)
  \right),
\end{equation}
where $\bs{C}$ is the correction term in \eqref{eq:correction.term}. This formulation involves both $\bs{r}_i$ and $\nabla_{\bs{\phi}}\bs{g}$ at $\bs{\phi}=\bs{0}$, that is, no sampling of the nuisance parameters is needed.

We follow \cite{Bec18} and estimate the EIG using the DLMC estimator

\begin{equation}
  \hat{I}_{\mathrm{DL}} = \frac{1}{N}\sum_{n=1}^N\left[\log\left(\tilde{p}(\bs{Y}^{(n)}|\bs{\theta}^{(n)})\right)-\log\left(\frac{1}{M}\sum_{m=1}^{M}\tilde{p}(\bs{Y}^{(n)}|\bs{\vartheta}^{(n,m)})\right)\right],
\end{equation}
where $\bs{\theta}^{(n)}\stackrel{\mathrm{iid}}{\sim}\pi(\bs{\theta})$ and  $\bs{Y}^{(n)}\stackrel{\mathrm{iid}}{\sim}\tilde{p}(\bs{Y}|\bs{\theta}^{(n)})$ for the first term and $\bs{\vartheta}^{(n,m)}\stackrel{\mathrm{iid}}{\sim}\pi(\bs{\theta})$ for the second term. We introduce this estimator only as a reference as the following estimators based on the Laplace approximation are much more computationally efficient and accurate.
When the forward model $\bs{g}$ includes a PDE, we use a discretized version $\bs{g}_h$ instead with mesh discretization parameter $h$. As $h \to 0$ asymptotically, the convergence order of $\bs{g}_h$ is given by
\begin{equation}
\bb{E}[\lVert \bs{g}(\bs{\theta})-\bs{g}_h(\bs{\theta})\rVert_2] = \cl{O}(h^\eta),
\end{equation}
where $\eta>0$ is the $h$-convergence rate in the weak sense. The work of evaluating $\bs{g}_h$ is assumed to be $\cl{O}(h^{-\gamma})$, for some $\gamma > 0$. We also assume that $\bs{g}$ is twice differentiable with respect to $\bs{\theta}$ and differentiable with respect to $\bs{\phi}$.

\subsection{Optimal setting for the Double-loop Monte Carlo estimator}
The average computational work of the DLMC estimator is
\begin{equation}\label{eq:Work}
  W_{\mathrm{DL}}\propto NMh^{-\gamma}.
\end{equation}
To derive the optimal setting for this estimator, we follow \cite{Bec18} and split the error into bias and variance, estimating each individually.
\begin{equation}
|\hat{I}_{\mathrm{DL}}-I|\leq\underbrace{|\mathbb{E}[\hat{I}_{\mathrm{DL}}]-I|}_{\text{bias error}}+\underbrace{|\hat{I}_{\mathrm{DL}}-\mathbb{E}[\hat{I}_{\mathrm{DL}}]|}_{\text{variance error}},
\end{equation}
\begin{equation}
  |\mathbb{E}[\hat{I}_{\mathrm{DL}}]-I|\leq C_{\mathrm{DL},3}h^{\eta}+\frac{C_{\mathrm{DL},4}}{M}+o(h^{\eta})+\cl{O}\left(\frac{1}{M^2}\right),
\end{equation}
\begin{equation}
  \bb{V}[\hat{I}_{\mathrm{DL}}]\leq \frac{C_{\mathrm{DL},1}}{N}+\frac{C_{\mathrm{DL},2}}{NM}+\cl{O}\left(\frac{1}{NM^2}\right),
\end{equation}
where
\begin{equation}
  C_{\mathrm{DL},1}=\bb{V}\left[\log\left(\frac{\tilde{p}(\bs{Y|\bs{\theta}})}{\tilde{p}(\bs{Y})}\right)\right],
\end{equation}
  \begin{equation}
    C_{\mathrm{DL},2}=\bb{E}\left[\bb{V}\left[\frac{\tilde{p}(\bs{Y|\bs{\theta}})}{\tilde{p}(\bs{Y})}|\bs{Y}\right]\right],
\end{equation}
where $C_{\mathrm{DL},3}$ is the constant of the $h$-convergence of $\hat{I}_{\mathrm{DL}}$ and
\begin{equation}
  C_{\mathrm{DL},4}=\frac{1}{2}\bb{E}\left[\bb{V}\left[\frac{\tilde{p}(\bs{Y|\bs{\theta}})}{\tilde{p}(\bs{Y})}|\bs{Y}\right]\right].
\end{equation}
We obtain the optimal setting by solving
\begin{equation}
  (N^{\ast},M^{\ast},h^{\ast},\kappa^{\ast})=\argmin_{(N,M,h,\kappa)}W_{\mathrm{DL}} \quad \text{subject to}\quad \begin{cases}\frac{C_{\mathrm{DL},1}}{N}+\frac{C_{\mathrm{DL},2}}{NM}\leq(\kappa TOL/C_{\alpha})^2\\C_{\mathrm{DL},3}h^{\eta}+\frac{C_{\mathrm{DL},4}}{M}\leq(1-\kappa)TOL\end{cases},
\end{equation}
for $C_{\alpha}=\bs{\Phi}(1-\frac{\alpha}{2})$, the inverse cumulative distribution function of the standard normal at confidence level $1-\alpha$ and $TOL>0$ is the allotted tolerance, which results in
\begin{equation}
  N^{\ast}=\frac{C_{\alpha}^2}{2\kappa^{\ast}}\frac{C_{\mathrm{DL},1}}{1-\kappa^{\ast}(1+\frac{\gamma}{2\eta})}TOL^{-2},
\end{equation}
\begin{equation}
  M^{\ast}=\frac{C_{\mathrm{DL},2}}{2\left(1-\kappa^{\ast}(1+\frac{\gamma}{2\eta})\right)}TOL^{-1},
\end{equation}
\begin{equation}
  h^{\ast}=\left(\frac{\gamma}{\eta}\frac{\kappa^{\ast}}{2C_{\mathrm{DL},3}}\right)^{1/\eta}TOL^{1/\eta}
\end{equation}
and the splitting parameter $\kappa^{\ast}\in (0,1)$ is given as the solution of the quadratic equation:
\begin{equation}
  \left[\frac{1}{C_{\mathrm{DL},1}}\left(1+\frac{\gamma}{2\eta}\right)^2TOL\right]\kappa^{\ast 2}-\left[\frac{1}{4}+\left(\frac{1}{2}+\frac{2}{C_{\mathrm{DL},1}}TOL\right)\left(1+\frac{\gamma}{2\eta}\right)\right]\kappa^{\ast}+\left[\frac{1}{2}+\frac{1}{C_{\mathrm{DL},1}}TOL\right]=0.
\end{equation}
The average work in the optimal setting is given as
\begin{equation}
  W_{\mathrm{DL}}^{\ast}\propto TOL^{-(3+\frac{\gamma}{\eta})}.
\end{equation}

\section{Monte Carlo Laplace approximation}\label{sec:laplace}

The DLMC estimator is expensive and often suffers from numerical underflow \cite{Bec18}. Some inexpensive and accurate estimators have been proposed for the EIG to reduce the computational cost: importance sampling \cite{Bec18} or multilevel techniques combined with importance sampling \cite{Bec20}. Furthermore, the inner integral may be estimated with the direct use of the Laplace approximation, tracking the correction term $\bs{C}(\bs{\theta})$ introduced by the small-noise approximation throughout.

\subsection{Monte Carlo Laplace estimator}
We approximate the posterior pdf of $\bs{\theta}$ using \eqref{eq:likelihood.tilde} as
\begin{align}\label{eq:posterior.small}
\pi(\bs{\theta}|\bs{Y})\approx&\frac{\tilde{p}(\bs{Y}|\bs{\theta})\pi(\bs{\theta})}{p(\bs{Y})}\nonumber\\
\propto &\det\left(2\pi\bs{\Sigma}_{\bs\epsilon}^{1/2}\bs{C}(\bs{\theta})^{-1}\bs{\Sigma}_{\bs\epsilon}^{1/2}\right)^{-\frac{N_e}{2}}\exp\left(-\frac{1}{2}\sum_{i=1}^{N_e}\bs{r}_i(\bs{\theta},\bs{0},\bs{\epsilon}_i)\cdot\bs{\Sigma}_{\bs\epsilon}^{-1/2}\bs{C}(\bs{\theta})\bs{\Sigma}_{\bs\epsilon}^{-1/2}\bs{r}_i(\bs{\theta},\bs{0},\bs{\epsilon}_i)
\right)\pi(\bs{\theta}).
\end{align}
We obtain a Gaussian approximation of this posterior from the Laplace approximation and consider the negative logarithm of \eqref{eq:posterior.small}
\begin{equation}
  F(\bs{\theta})\coloneqq \frac{1}{2}\sum_{i=1}^{N_e}\bs{r}_i(\bs{\theta},\bs{0},\bs{\epsilon}_i)\cdot\bs{\Sigma}_{\bs\epsilon}^{-1/2}\bs{C}(\bs{\theta})\bs{\Sigma}_{\bs\epsilon}^{-1/2}\bs{r}_i(\bs{\theta},\bs{0},\bs{\epsilon}_i)-h(\bs{\theta})+D(\bs{\theta}),
\end{equation}
where $h(\bs{\theta})=\log(\pi(\bs{\theta}))$ and
\begin{equation}
    D(\bs{\theta})=\log\left(\det\left(2\pi\bs{\Sigma}_{\bs\epsilon}^{1/2}\bs{C}(\bs{\theta})^{-1}\bs{\Sigma}_{\bs\epsilon}^{1/2}\right)^{\frac{N_e}{2}}\right).
   \end{equation}
Next, we take the second-order Taylor approximation of $F$ around the maximum a posteriori (MAP) estimate $\hat{\bs{\theta}}$
\begin{equation}\label{eq:taylor}
  \tilde{F}(\bs{\theta})=F(\hat{\bs{\theta}})+\nabla_{\bs{\theta}}F(\hat{\bs{\theta}})(\bs{\theta}-\hat{\bs{\theta}})+\frac{1}{2}(\bs{\theta}-\hat{\bs{\theta}})\cdot\nabla_{\bs{\theta}}\nabla_{\bs{\theta}}F(\hat{\bs{\theta}})(\bs{\theta}-\hat{\bs{\theta}}),
\end{equation}
where
\begin{multline}
\nabla_{\bs{\theta}}F(\hat{\bs{\theta}})=-\nabla_{\bs{\theta}}\bs{g}(\hat{\bs{\theta}},\bs{0})\cdot\bs{\Sigma}_{\bs\epsilon}^{-1/2}\bs{C}(\hat{\bs{\theta}})\bs{\Sigma}_{\bs\epsilon}^{-1/2}\sum_{i=1}^{N_e}\bs{r}_i(\hat{\bs{\theta}},\bs{0},\bs{\epsilon}_i)\\
+ \frac{1}{2}\sum_{i=1}^{N_e}\bs{r}_i(\hat{\bs{\theta}},\bs{0},\bs{\epsilon}_i)\cdot\bs{\Sigma}_{\bs\epsilon}^{-1/2}\nabla_{\bs{\theta}}\bs{C}(\hat{\bs{\theta}})\bs{\Sigma}_{\bs\epsilon}^{-1/2}\bs{r}_i(\hat{\bs{\theta}},\bs{0},\bs{\epsilon}_i)-\nabla_{\bs{\theta}}h(\hat{\bs{\theta}})+\nabla_{\bs{\theta}}D(\hat{\bs{\theta}})
\end{multline}
which is zero because it is evaluated at the MAP, and
\begin{align}\label{eq:sigma}
\nabla_{\bs{\theta}}\nabla_{\bs{\theta}}F(\hat{\bs{\theta}})={}&-\nabla_{\bs{\theta}}\nabla_{\bs{\theta}}\bs{g}(\hat{\bs{\theta}},\bs{0})\cdot\bs{\Sigma}_{\bs\epsilon}^{-1/2}\bs{C}(\hat{\bs{\theta}})\bs{\Sigma}_{\bs\epsilon}^{-1/2}\sum_{i=1}^{N_e}\bs{r}_i(\hat{\bs{\theta}},\bs{0},\bs{\epsilon}_i)\nonumber\\[4pt]
&-2\nabla_{\bs{\theta}}\bs{g}(\hat{\bs{\theta}},\bs{0})\cdot\bs{\Sigma}_{\bs\epsilon}^{-1/2}\nabla_{\bs{\theta}}\bs{C}(\hat{\bs{\theta}})\bs{\Sigma}_{\bs\epsilon}^{-1/2}\sum_{i=1}^{N_e}\bs{r}_i(\hat{\bs{\theta}},\bs{0},\bs{\epsilon}_i)\nonumber\\[4pt]
&+N_e\nabla_{\bs{\theta}}\bs{g}(\hat{\bs{\theta}},\bs{0})\cdot\bs{\Sigma}_{\bs\epsilon}^{-1/2}\bs{C}(\hat{\bs{\theta}})\bs{\Sigma}_{\bs\epsilon}^{-1/2}\nabla_{\bs{\theta}}\bs{g}(\hat{\bs{\theta}},\bs{0})\nonumber\\[4pt]
&+ \frac{1}{2}\sum_{i=1}^{N_e}\bs{r}_i(\hat{\bs{\theta}},\bs{0},\bs{\epsilon}_i)\cdot\bs{\Sigma}_{\bs\epsilon}^{-1/2}\nabla_{\bs{\theta}}\nabla_{\bs{\theta}}\bs{C}(\hat{\bs{\theta}})\bs{\Sigma}_{\bs\epsilon}^{-1/2}\bs{r}_i(\hat{\bs{\theta}},\bs{0},\bs{\epsilon}_i)-\nabla_{\bs{\theta}}\nabla_{\bs{\theta}}h(\hat{\bs{\theta}})+\nabla_{\bs{\theta}}\nabla_{\bs{\theta}}D(\hat{\bs{\theta}}).
\end{align}
The MAP estimate $\hat{\bs{\theta}}$ of \eqref{eq:posterior.small} is given by
\begin{align}\label{eq:theta.hat}
\hat{\bs{\theta}}&=\argmax_{\bs{\theta}}\pi(\bs{\theta}|\bs Y)\approx\argmin_{\bs{\theta}}F(\bs{\theta})\nonumber\\
&=\argmin_{\bs{\theta}}\left[\frac{1}{2}\sum_{i=1}^{N_e}\bs{r}_i(\bs{\theta},\bs{0},\bs{\epsilon}_i) \cdot \bs{\Sigma}_{\bs{\varepsilon}}^{-1/2}\bs{C}(\bs{\theta})\bs{\Sigma}_{\bs{\varepsilon}}^{-1/2}\bs{r}_i(\bs{\theta},\bs{0},\bs{\epsilon}_i)-h(\bs{\theta})-\frac{N_e}{2}\log\det(\bs{C}(\bs{\theta}))\right].
\end{align}
\begin{rmk}[$\hat{\bs{\theta}}-\bs{\theta}_t=\cl{O}_{\bb{P}}\left(\frac{1}{\sqrt{N_e}}\right)$]\label{rmk:foa}
From the first-order approximation of $\hat{\bs{\theta}}$ (Appendix \ref{ap:map}),
\begin{equation}\label{eq:theta.true}
\hat{\bs{\theta}}=\bs{\theta}_t+\cl{O}_{\bb{P}}\left(\frac{1}{\sqrt{N_e}}\right).
\end{equation}
Thus, $\hat{\bs{\theta}}$ can be approximated by $\bs{\theta}_t$ for small $\lVert\bs{\Sigma}_{\bs{\phi}}\rVert$ and large $N_e$.
\end{rmk}
Given Remark \ref{rmk:foa}, we obtain $\sum_{i=1}^{N_e}\bs{r}_i(\hat{\bs{\theta}},\bs{0},\bs{\epsilon}_i)\sim\cl{N}(N_e(\bs{g}(\bs{\theta}_t,\bs{0})-\bs{g}(\hat{\bs{\theta}},\bs{0})),N_e\bs{\Sigma}_{\bs{\epsilon}})=\cl{O}_{\bb{P}}(\sqrt{N_e})$, so we drop the first, second, and fourth terms in \eqref{eq:sigma}.

We obtain the Laplace approximation of the posterior as
\begin{equation}\label{eq:lap.post}
  \tilde{\pi}(\bs{\theta}|\bs{Y})=\frac{1}{(2\pi)^{\frac{d_{\bs{\theta}}}{2}}|\bs{\Sigma}|^{\frac{1}{2}}}\exp\left(-\frac{(\bs{\theta}-\bs{\hat\theta}) \cdot \bs{\Sigma}^{-1}(\bs{\theta}-\bs{\hat\theta})}{2}\right),
\end{equation}
where
\begin{equation}\label{eq:sigma.la}
  \bs{\Sigma}^{-1}=N_e\nabla_{\bs{\theta}}\bs{g}(\hat{\bs{\theta}},\bs{0})\cdot\bs{\Sigma}_{\bs\epsilon}^{-1/2}\bs{C}(\hat{\bs{\theta}})\bs{\Sigma}_{\bs\epsilon}^{-1/2}\nabla_{\bs{\theta}}\bs{g}(\hat{\bs{\theta}},\bs{0})-\nabla_{\bs{\theta}}\nabla_{\bs{\theta}}h(\hat{\bs{\theta}}) +\nabla_{\bs{\theta}}\nabla_{\bs{\theta}}D(\hat{\bs{\theta}})+\cl{O}(\sqrt{N_e}),
\end{equation}
by taking the negative exponential of \eqref{eq:taylor}. We follow \cite{Lon13} to derive a sample-based estimator for the EIG, where all information about the nuisance parameter is encoded in the covariance matrix $\bs{\Sigma}$. Using the Laplace approximation \eqref{eq:lap.post} to rewrite the log-ratio appearing in the EIG \eqref{eq:EIG.definition}, we obtain
\begin{align}
  \log\left(\frac{\pi(\bs{\theta}|\bs Y)}{\pi(\bs{\theta})}\right)&=\underbrace{\log\left(\frac{\pi(\bs{\theta}|\bs Y)}{\tilde{\pi}(\bs{\theta}|\bs Y)}\right)}_{:=\epsilon_{\mathrm{LA}}}+\log\left(\frac{\tilde{\pi}(\bs{\theta}|\bs Y)}{\pi(\bs{\theta})}\right)\nonumber\\
  &=\epsilon_{\mathrm{LA}}-\frac{1}{2}\log((2\pi)^{d_{\bs{\theta}}}|\bs{\Sigma}|)-\left(\frac{(\bs{\theta}-\bs{\hat\theta}) \cdot \bs{\Sigma}^{-1}(\bs{\theta}-\bs{\hat\theta})}{2}\right)-\log(\pi(\bs{\theta})).
\end{align}
From the Kullback--Leibler divergence \eqref{Dkl}, we obtain
\begin{align}
	D_{\mathrm{KL}}&=\int_{\bs\Theta}\log\left(\frac{\pi(\bs{\theta}|\bs Y)}{\pi(\bs{\theta})}\right)\tilde{\pi}(\bs{\theta}|\bs Y)\di{}\bs\theta+\underbrace{\int_{\bs\Theta}\log\left(\frac{\pi(\bs{\theta}|\bs Y)}{\pi(\bs{\theta})}\right)(\pi(\bs{\theta}|\bs Y)-\tilde{\pi}(\bs{\theta}|\bs Y))\di{}\bs\theta}_{\coloneqq \epsilon_{int}}\nonumber\\
	&=\int_{\bs\Theta}\epsilon_{\mathrm{LA}}\tilde{\pi}(\bs{\theta}|\bs Y)\di{}\bs\theta+\int_{\bs\Theta}\left[-\frac{1}{2}\log((2\pi)^{d_{\bs{\theta}}}|\bs{\Sigma}|)-\left(\frac{(\bs{\theta}-\bs{\hat\theta}) \cdot \bs{\Sigma}^{-1}(\bs{\theta}-\bs{\hat\theta})}{2}\right)-\log(\pi(\bs{\theta}))\right]\tilde{\pi}(\bs{\theta}|\bs Y)\di{}\bs\theta+\epsilon_{int}\nonumber\\
	&=-\frac{1}{2}\log((2\pi)^{d_{\bs{\theta}}}|\bs{\Sigma}|)-\frac{d_{\bs\theta}}{2}-\log(\pi(\hat{\bs\theta}))-\frac{\text{tr}(\bs\Sigma\nabla_{\bs\theta}\nabla_{\bs\theta}\log(\pi(\hat{\bs\theta})))}{2}+\cl{O}_{\bb{P}}\left(\frac{1}{N_e^2}\right),
\end{align}
where we used
\begin{equation}
 \int_{\bs\Theta}\log(\pi(\bs{\theta}))\tilde{\pi}(\bs{\theta}|\bs Y)\di{}\bs\theta=\log(\pi(\hat{\bs\theta}))+\frac{\text{tr}(\bs\Sigma\nabla_{\bs\theta}\nabla_{\bs\theta}\log(\pi(\hat{\bs\theta})))}{2}+\cl{O}_{\bb{P}}\left(\frac{1}{N_e^2}\right),
\end{equation}
\begin{equation}
  \int_{\bs\Theta}\epsilon_{\mathrm{LA}}\tilde{\pi}(\bs{\theta}|\bs Y)\di{}\bs\theta = \cl{O}_{\bb{P}}\left(\frac{1}{N_e^2}\right),
\end{equation}
and
\begin{equation}
   \epsilon_{int}=\cl{O}_{\bb{P}}\left(\frac{1}{N_e^2}\right).
\end{equation}
For the detailed derivation, see Appendices A--C of \cite{Lon13}. The difference between the estimator constructed in \cite{Lon13} and ours lies in the covariance matrix $\bs{\Sigma}$. Our covariance matrix also contains information on the nuisance uncertainty. However, in both cases, the covariance matrices converge to the true one at the same rate in probability $\cl{O}_{\bb{P}}\left(\frac{1}{N_e^2}\right)$.

Taking the expected value over all $\bs Y$ and using the law of total probability results in the EIG
\begin{equation}
	I=\int_{\bs\Theta}\int_{\cl{Y}}\left[-\frac{1}{2}\log((2\pi)^{d_{\bs{\theta}}}|\bs{\Sigma}|)-\frac{d_{\bs\theta}}{2}-\log(\pi(\hat{\bs\theta}))-\frac{\text{tr}(\bs\Sigma\nabla_{\bs\theta}\nabla_{\bs\theta}\log(\pi(\hat{\bs\theta})))}{2}\right]p(\bs Y|\bs\theta_t)\di{}\bs Y\pi(\bs\theta_t)\di{}\bs\theta_t+\cl{O}_{\bb{P}}\left(\frac{1}{N_e^2}\right).
\end{equation}
Using \eqref{eq:theta.true}, we obtain the sample-based estimator
\begin{equation}\label{eq:mcla_sample}
  \hat{I}_{\mathrm{LA}}=\frac{1}{N}\sum_{n=1}^N\left(-\frac{1}{2}\log((2\pi)^{d_{\bs{\theta}}}|\bs{\Sigma}(\bs{\theta}^{(n)})|)-\frac{d_{\bs{\theta}}}{2}-\log(\pi(\bs\theta^{(n)}))\right),
\end{equation}
where $\bs\theta^{(n)}\stackrel{\mathrm{iid}}{\sim} \pi(\bs\theta)$ are sampled from the prior. This estimator has a bias error of $\cl{O}_{\bb{P}}\left(\frac{1}{N_e}\right)$.

\subsection{Optimal setting for the Monte Carlo Laplace estimator}
The average computational work of the estimator is
\begin{equation}
W_{\mathrm{LA}}\propto NJh^{-\gamma},
\end{equation}
where $J$ is the number of forward model evaluations required for the Jacobian matrix in \eqref{eq:theta.hat} and \eqref{eq:sigma.la}.

To obtain the optimal setting, we follow \cite{Bec18}, splitting the error into the bias and variance:
\begin{equation}
  |\mathbb{E}[\hat{I}_{\mathrm{LA}}]-I|\leq C_{\mathrm{LA},3}h^{\eta}+\frac{C_{\mathrm{LA},2}}{N_e}+o(h^{\eta}),
\end{equation}
\begin{equation}
  \bb{V}[\hat{I}_{\mathrm{LA}}]\leq \frac{C_{\mathrm{LA},1}}{N},
\end{equation}
where
\begin{equation}
  C_{\mathrm{LA},1}=\bb{V}[D_{\mathrm{KL}}].
\end{equation}

Then, we solve
\begin{equation}
  (N^{\ast},h^{\ast},\kappa^{\ast})=\argmin_{(N,h,\kappa)}W_{\mathrm{LA}} \quad \text{subject to}\quad \begin{cases}\frac{C_{\mathrm{LA},1}}{N}\leq(\kappa TOL/C_{\alpha})^2\\C_{\mathrm{LA},3}h^{\eta}+\frac{C_{\mathrm{LA},2}}{N_e}\leq(1-\kappa)TOL\end{cases},
\end{equation}
which is only solvable if the bias introduced by the Laplace approximation does not exceed the tolerance $TOL$, that is,
\begin{equation}\label{eq:laplace.bias}
  \frac{C_{\mathrm{LA},2}}{N_e}\leq(1-\kappa)TOL.
\end{equation}
If this condition is met, the solution is given by
\begin{equation}
N^{\ast}=C_{\rm{LA},1}\left(\frac{C_{\alpha}}{\kappa^{\ast}}\right)^2TOL^{-2},
\end{equation}
\begin{equation}\label{eq:mcla.h}
  h^{\ast}=\left(\frac{\gamma}{\eta}\frac{\kappa^{\ast}}{2C_{\mathrm{LA},3}}\right)^{1/\eta}TOL^{1/\eta},
\end{equation}
\begin{equation}
  \kappa^{\ast}=\frac{1-\frac{C_{\mathrm{LA},2}}{N_e TOL}}{1+\frac{\gamma}{2\eta}}.
\end{equation}
The average work in the optimal setting is given as
\begin{equation}
  W_{\mathrm{LA}}^{\ast}\propto TOL^{-(2+\frac{\gamma}{\eta})};
\end{equation}
thus, a reduction in the average work is achieved compared to the DLMC estimator, but a lower bound for the allowed tolerance is introduced.

\section{Laplace-based importance sampling approximation}\label{sec:importance.sampling}

Another method of estimating the EIG \eqref{eq:EIG.definition} is to use importance sampling, which amounts to a change of measure for Monte Carlo sampling in the inner loop. As the new measure, we use \eqref{eq:lap.post} conditioned on the likelihood $\bs{Y}$ obtained from the outer loop. From the small-noise approximation \eqref{eq:likelihood.tilde}, we obtain
\begin{equation}
  I=\int_\Theta\int_\cl{Y}\left[\log(\tilde{p}(\bs{Y}|\bs{\theta}))-\log\left(\int_\Theta \frac{\tilde{p}(\bs{Y}|\bs{\vartheta})\pi(\bs{\vartheta})}{\tilde{\pi}(\bs{\vartheta}|\bs{Y})}\tilde{\pi}(\bs{\vartheta}|\bs{Y})\di{}\bs{\vartheta}\right)\right]\tilde{p}(\bs{Y}|\bs{\theta})\di{}\bs{Y}\pi(\bs{\theta})\di{}\bs{\theta}
\end{equation}
with the sample-based estimator
\begin{equation}
  \hat{I}_{\mathrm{DLIS}}=\frac{1}{N}\sum_{n=1}^N\left[\log\left(
\tilde{p}(\bs{Y}^{(n)}|\bs{\theta}^{(n)})\right)-\log\left(
\frac{1}{M}\sum_{m=1}^M\frac{
\tilde{p}(\bs{Y}^{(n)}|\bs{\vartheta}^{(n,m)})\pi(\bs{\vartheta}^{(n,m)})}
{\tilde{\pi}(\bs{\vartheta}^{(n,m)}|\bs{Y}^{(n)})}
\right)\right],
\end{equation}
where $\bs{\theta}^{(n)}\stackrel{\mathrm{iid}}{\sim}\pi(\bs{\theta})$ and  $\bs{Y}^{(n)}\stackrel{\mathrm{iid}}{\sim}\tilde{p}(\bs{Y}|\bs{\theta}^{(n)})$ for the first term and $\bs{\vartheta}^{(n,m)}\stackrel{\mathrm{iid}}{\sim}\tilde{\pi}(\bs{\theta}|\bs{Y}^{(n)})$ for the second term.

\begin{rmk}[Optimal setting of the DLMCIS estimator]
The derivation for the optimal work of the DLMCIS estimator is the same, if we neglect the cost of computing the new measure, as for the case without importance sampling. The only difference is that the constants depending on the variance are significantly smaller, decreasing the samples needed, especially for the optimal number of inner samples $M^{\ast}$. The bias in this method is controllable by the number of inner samples, in contrast to the MCLA estimator.
\end{rmk}

\section{Numerical results}\label{sec:numerical.results}

\subsection{Example 1: Small-noise approximation degeneration}

\begin{figure}[ht]
    \centering
    \includegraphics[width=0.6\textwidth]{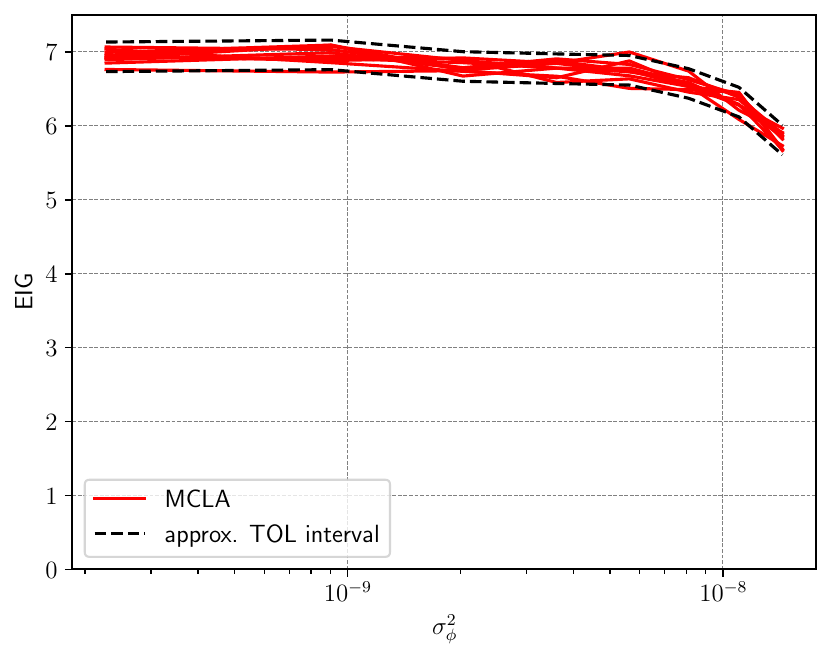}
    \caption{Example 1: EIG vs. $\sigma_{\phi}^2$ for 10 runs of the MCLA estimator with approximate tolerance interval. The covariance matrix in \eqref{eq:small.likelihood.fam} loses positive definiteness when $\sigma_{\phi}^2$ exceeds a certain threshold of approximately $1.44\times 10^{-8}$.}
    \label{fig:breakdown}
\end{figure}

To depict the range of applicability of small-noise approximation, we studied the following model with respect to the norm of the covariance of the nuisance parameters:
\begin{equation}
	\bs{g}(\bs{\theta},\bs{\phi})=\bs{A}\bs{\theta}(1+\bs{\phi}\cdot \bs{e}),
\end{equation}
where
\begin{equation}
	\bs{A}=\begin{pmatrix}
10 & 2 \\
5 & 20 \\
50 & -2
\end{pmatrix}, \quad \bs{e}=\begin{pmatrix}
	1 \\ 1
\end{pmatrix},
\end{equation}
and $\bs{\epsilon}\sim\cl{N}(\bs{0}, \sigma_{\epsilon}^2\bs{1})$, $\bs{\theta}\sim\cl{N}(\bs{e}, \sigma_{\theta}^2\bs{1})$, $\bs{\phi}\sim\cl{N}(\bs{0}, \sigma_{\phi}^2\bs{1})$, $\sigma_{\epsilon}^2=\sigma_{\theta}^2=10^{-4}$ and $\bs{1}$ is the identity matrix in $\bb{R}^2$. When varying $\sigma_{\phi}^2$, the MCLA estimator cannot handle values larger than $1.44\times 10^{-8}$. When the variance of the nuisance uncertainty $\sigma_{\phi}^2$ becomes larger than this threshold value, the correction term $\bs{C}$ in \eqref{eq:likelihood.tilde} causes the covariance matrix in \eqref{eq:small.likelihood.fam} to lose positive definiteness (see Remark \ref{rmk.eig}), yielding negative infinity in the logarithm expression in \eqref{eq:mcla_sample}. This loss results in the EIG approximation deteriorating and eventually breaking down completely. The same also happens for the DLMC and DLMCIS estimators. Figure \ref{fig:breakdown} illustrates 10 runs of the MCLA estimator with error tolerance specified as $TOL = 0.2$. Also, note that we do not account for the bias resulting from \eqref{eq:small.exp} for our numerical examples, since its computation is very demanding.

\begin{figure}[ht]
	\subfloat[Small $\sigma_{\phi}^2$\label{subfig-1:comp.small}]{%
		\includegraphics[width=0.45\textwidth]{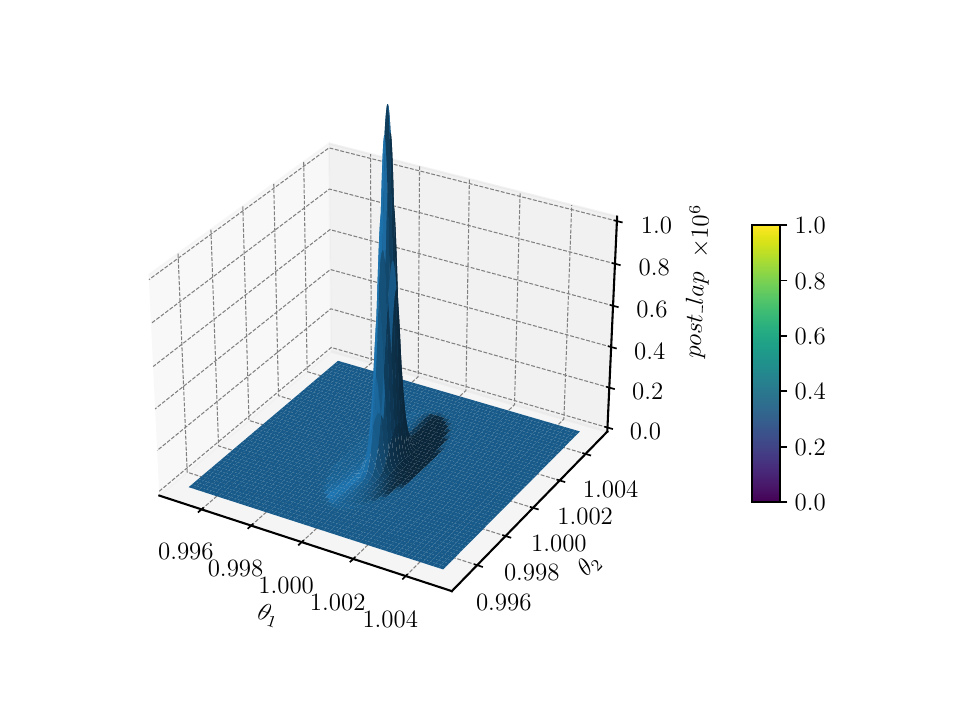}
	}
	\subfloat[Large $\sigma_{\phi}^2$\label{subfig-2:comp.large}]{%
		\includegraphics[width=0.45\textwidth]{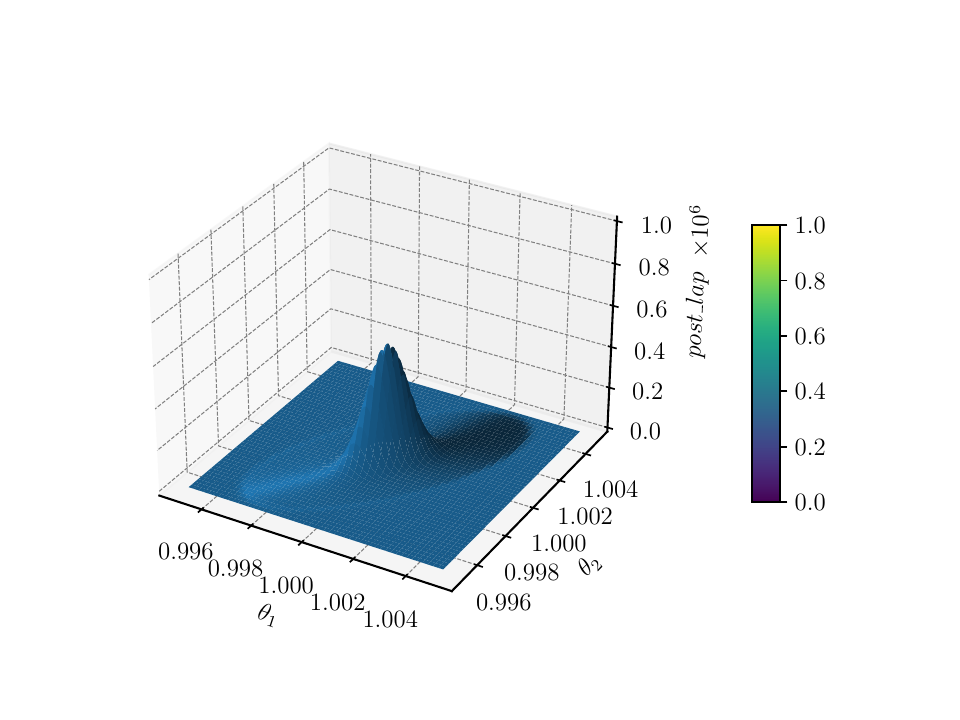}
	}
	\caption{Example 1: Laplace approximation of the posterior pdf is more concentrated for small nuisance uncertainty $\sigma_{\phi}^2=1.44\times 10^{-16}$ (left) than for a larger nuisance uncertainty $\sigma_{\phi}^2=1.44\times 10^{-8}$ (right).}
	\label{fig:eit.post.comp}
\end{figure}

The posterior covariance matrix computed from the MCLA estimator depends on the magnitude of the nuisance uncertainty $\sigma_{\phi}^2$ as indicated in Figure \ref{fig:eit.post.comp}. We present the Gaussian posterior pdf once for small $\sigma_{\phi}^2=1.44\times 10^{-16}$ and once for larger $\sigma_{\phi}^2=1.44\times 10^{-8}$ where the small-noise approximation starts breaking down. Naturally, we obtain a more concentrated posterior for the smaller nuisance uncertainty, meaning that the experiment is more informative in the absence of nuisance uncertainty.

\subsection{Example 2: Linear Gaussian example}

\begin{figure}[ht]
\begin{tikzpicture}
\node[inner sep=0pt] (comp) at (0,0)
    {\includegraphics[width=0.6\textwidth]{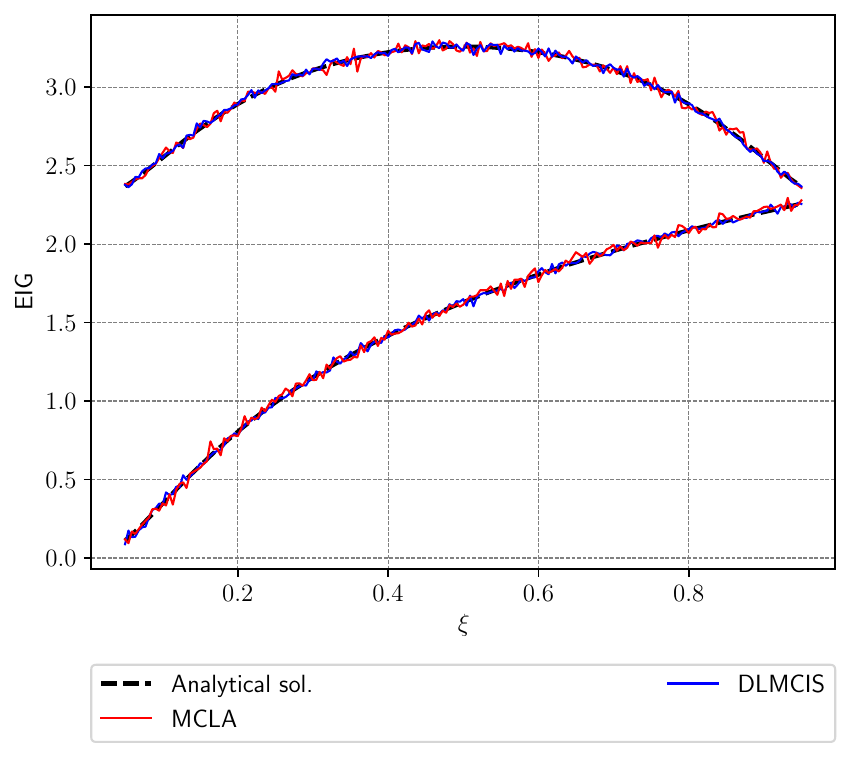}};
		\node at (0.65,2.6) {\small No nuisance uncertainty in the EIG};
    \node at (0.65,-1.0) {\small With nuisance uncertainty in the EIG};
\end{tikzpicture}
\caption{Example 2: EIG vs. design parameter $\xi$. Analytical solution (dashed black), DLMCIS estimator (solid blue), and MCLA estimator (solid red) for the case with or without nuisance uncertainty.}
\label{fig:fm_2}
\end{figure}

The following problem is constructed based on an example by \cite{Fen19}, and we use it to assess the efficiency of different estimators because it can be solved analytically. The linear Gaussian example has the form
\begin{equation}
	\bs{g}(\theta,\phi)=\bs{A}\bs{z},
\end{equation}
where
\begin{equation}
	\bs{A}=\begin{pmatrix}
\xi & 0 \\
 0 & (1-\xi)
\end{pmatrix}, \quad \bs{z}=\begin{pmatrix}
	\theta \\ \phi
\end{pmatrix},
\end{equation}
and $\bs{\epsilon}\sim\cl{N}(\bs{0}, \sigma_{\epsilon}^2\bs{1})$, $\theta\sim\cl{N}(0, \sigma_{\theta}^2)$, $\phi\sim\cl{N}(0, \sigma_{\phi}^2)$, $\sigma_{\epsilon}^2=\sigma_{\phi}^2=10^{-2}$, $\sigma_{\theta}^2=1$ and $\bs{1}$ is the identity matrix in $\bb{R}^2$. Alternatively, $\bs{\theta}$ can be chosen from $\bb{R}^2$ to observe how the optimal design shifts if there are no nuisance parameters. In this case, $\bs{z}=\bs{\theta}\sim\cl{N}(\bs{0},\bs{1})$. Figure \ref{fig:fm_2} indicates that both the DLMCIS and MCLA estimator accurately approximate the analytical solution for the case with nuisance parameters and without as long as $\sigma_{\phi}^2$ is kept small enough. As mentioned, if $\sigma_{\phi}^2$ is large, the correction term in the small-noise approximation causes the covariance matrix to lose positive definiteness. In this example, the model $\bs{g}$ is not bijective for $\xi=0$ and $\xi=1$. With these design choices, no inference can be made about $\bs{\theta}$. The DLMCIS estimator is less robust than the MCLA estimator close to these points.

\subsubsection{Optimal number of samples and error convergence}

\begin{figure}[ht]
	\subfloat[With nuisance uncertainty\label{subfig-1:opt.n}]{%
		\includegraphics[width=0.45\textwidth]{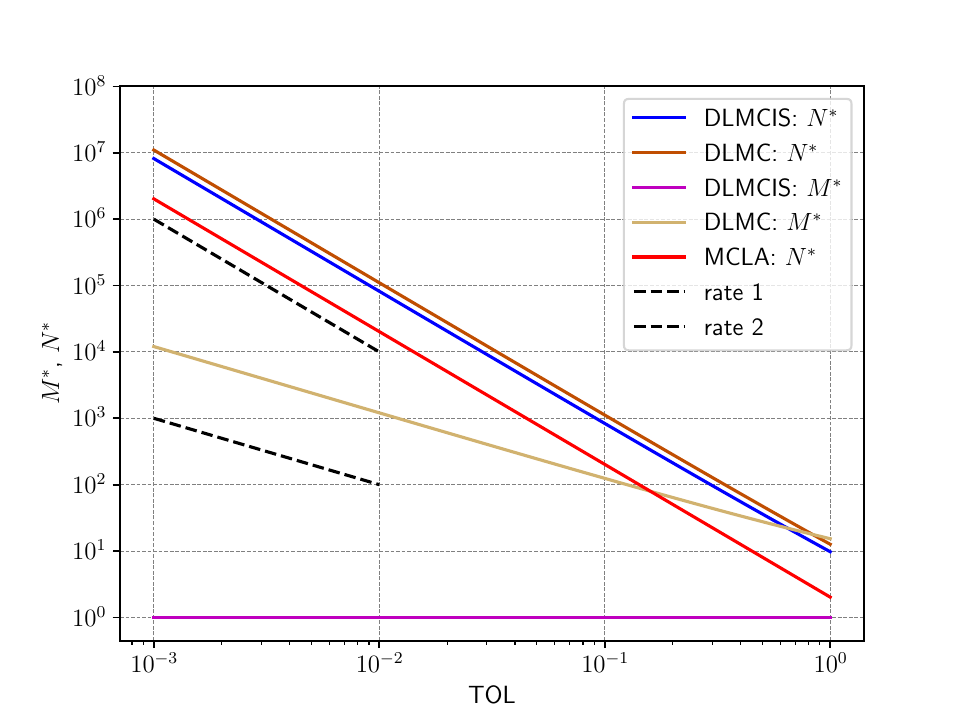}
	}
	\hfill
	\subfloat[No nuisance uncertainty\label{subfig-2:opt.nn}]{%
		\includegraphics[width=0.45\textwidth]{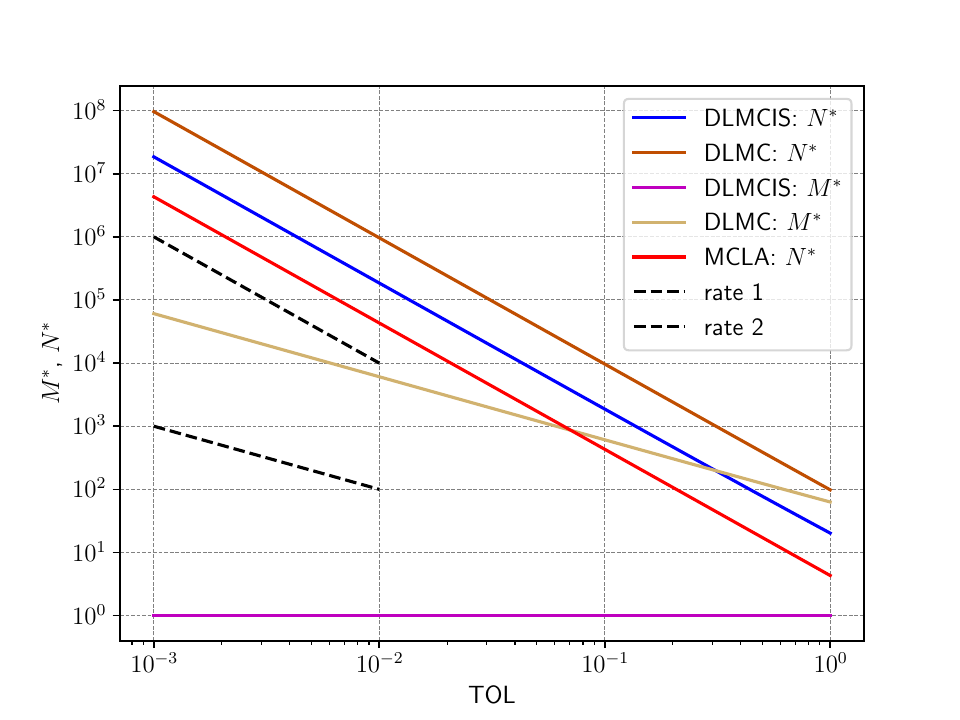}
	}
	\caption{Example 2: Optimal number of outer ($N^{\ast}$) and inner ($M^{\ast}$) samples vs. tolerance $TOL$. $M^{\ast}$ for the DLMCIS estimator starts increasing for small tolerances.}
	\label{fig:fm.opt}
\end{figure}

Figure \ref{fig:fm.opt} depicts the optimal number of outer samples $N^{\ast}$ and inner samples $M^{\ast}$ when computing the EIG for a fixed $\xi=1/2$. For the DLMC and DLMCIS estimator, the constants $C_{\mathrm{DL},1}$ and $C_{\mathrm{DL},2}$ were estimated using DLMC and DLMCIS, respectively, with $N=1000$ and $M=200$. For the MCLA estimator, the constant $C_{\mathrm{LA},1}$ was estimated using MCLA with $N=1000$.

 The optimal number of inner samples $M^{\ast}$ for DLMCIS increases at a rate of $\cl{O}(\text{TOL}^{-1})$ for small tolerances and is equal to $1$ for TOL greater than approximately $10^{-1}$. The optimal number of outer samples $N^{\ast}$ for all three estimators increases at a rate of $\cl{O}(\text{TOL}^{-2})$. This holds true for the case with and without nuisance uncertainty. Thus, the total complexity remains unchanged when including nuisance uncertainty into the model.

\begin{figure}[ht]
	\subfloat[DLMCIS\label{subfig-1:evt.dlis}]{%
		\includegraphics[width=0.45\textwidth]{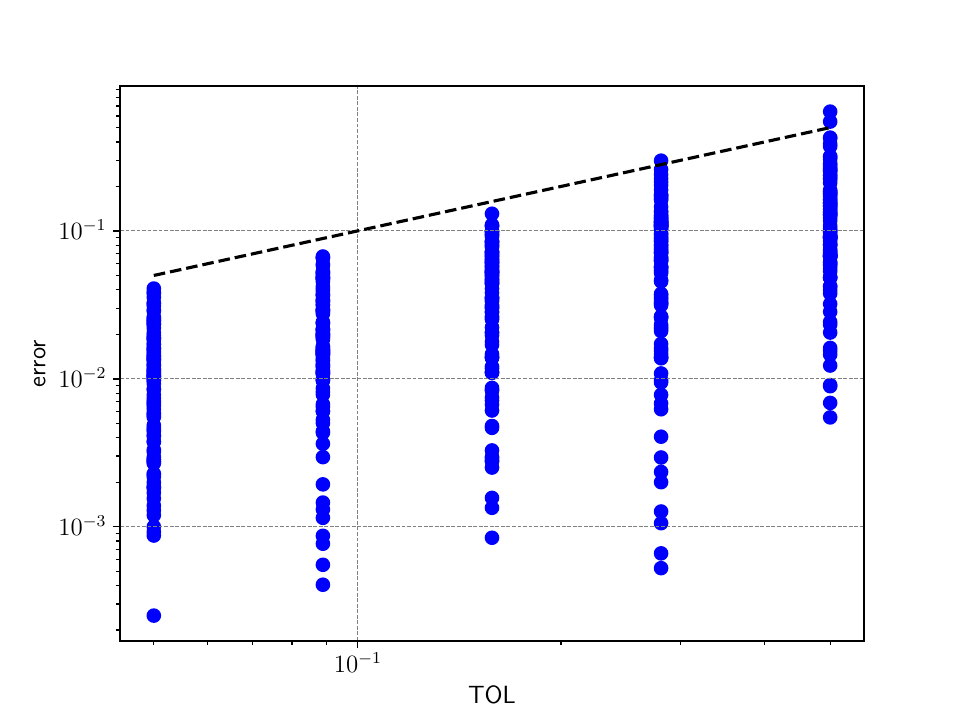}
	}
	\hfill
	\subfloat[MCLA\label{subfig-2:evt.mcla}]{%
		\includegraphics[width=0.45\textwidth]{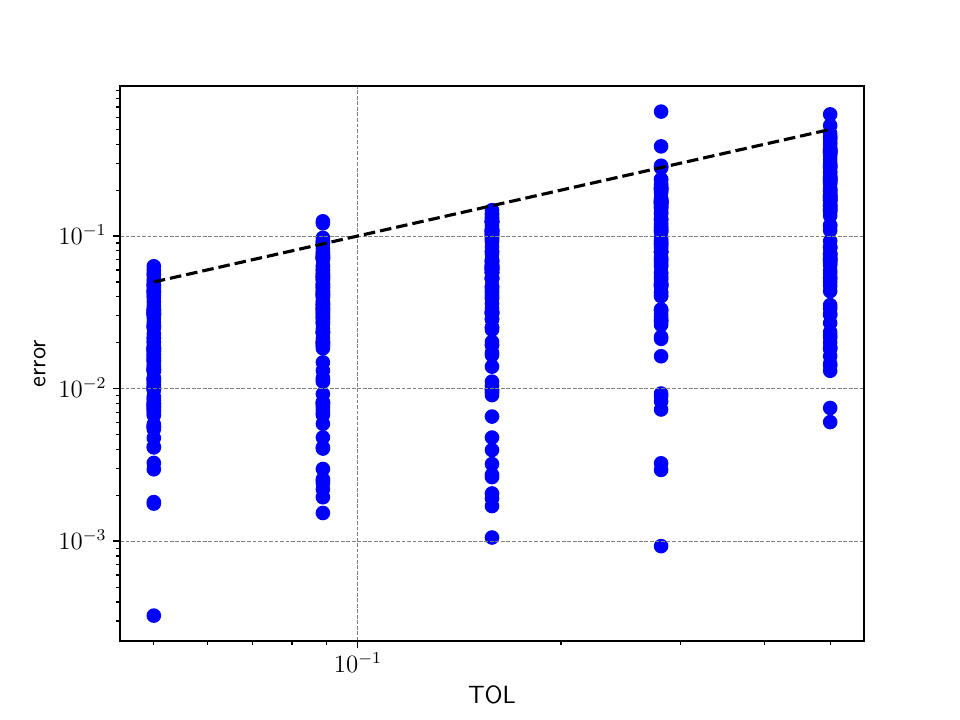}
	}
	\caption{Example 2: Error vs. tolerance consistency plot. Both estimators stay within the allowed tolerance $TOL$ with a predefined confidence parameter $C_\alpha=1.96$ ($\alpha=0.05$).}
	\label{fig:evt}
\end{figure}

The absolute error between the true solution and those estimated using MCLA and DLMCIS is presented in Figure \ref{fig:evt} for the case with nuisance uncertainty. The previously computed optimal number of samples was used. The consistency test indicates that the absolute error is lower than TOL with $95\%$ probability.

\subsection{Example 3: Electrical impedance tomography: Recovering fiber orientation in composite laminate materials under nuisance uncertainty}\label{sec:EIT}
We apply the proposed estimators to EIT, where low-frequency electric currents are injected through electrodes attached to the surface of a closed body to recover the fiber orientation in each orthotropic ply of a composite laminate material based on potential measurements. We adopt the complete electrode model from \cite{Som92}. The potential field is modeled by a second-order PDE with an electrode boundary model. We assume the exact interior conductivity to be unknown and it is therefore modeled using nuisance parameters. The goal is to determine the design of the experiment, namely, the shift and distance between electrodes, that allows for optimal recovery of the fiber angles under challenging conditions.

\subsubsection{EIT formulation}

The rectangular body $\cl{D}$ with boundary $\partial \cl{D}$ consists of two plies $\cl{D}=\bigcup_{i=1}^2\cl{D}_i$. The surfaces of the exterior boundary are equipped with $N_{\text{el}}$ electrodes $E_l$, $1\leq l\leq N_{\text{el}}$ with surface impedance $z_l$. We have the following PDE formulation
\begin{equation}\label{eq:pde.1}
	\nabla\cdot \bs{j}(\bs{x}, \omega)=0\quad \text{in}\, \cl{D}\, \text{and}
\end{equation}
\begin{equation}\label{eq:pde.2}
\bs{j}(\bs{x}, \omega)=	\bs{\bar{\sigma}}\cdot\nabla u(\bs{x}, \omega),
\end{equation}
where $u$ is the quasi-static potential field, $\bs{j}$ indicates the current flux, and $\bs{\bar{\sigma}}$ denotes the conductivity field given by
\begin{equation}
\bs{\bar{\sigma}}(\bs{x}, \omega)=\bs{Q}(\theta_i(\omega))^{\trans}\cdot \bs{\sigma}(\phi_i(\omega))\cdot\bs{Q}(\theta_i(\omega)), \quad \text{for}\, \bs{x}\in\cl{D}_i, \, i=1,2,
\end{equation}
where
\begin{equation}
  \bs{Q}(\theta_i) = \begin{pmatrix}
    \cos(\theta_i)& 0 & -\sin(\theta_i) \\ 0 & 1 & 0 \\ \sin(\theta_i) & 0 & \cos(\theta_i)
  \end{pmatrix} \quad \text{and} \quad \bs{\sigma}(\phi_i)=\begin{pmatrix}
    \sigma_1(\phi_i)& 0 & 0 \\ 0 & \sigma_2(\phi_i) & 0 \\ 0 & 0 & \sigma_3(\phi_i)
  \end{pmatrix},\, i=1,2.
\end{equation}
The boundary conditions are given by
\begin{equation}
	\begin{cases}\label{eq:pde.bc.1}
		\bs{j}\cdot\bs{n}=0 & \text{on}\, \partial\cl{D}\setminus (\cup E_l),\\
		\int_{E_l}\bs{j}\cdot\bs{n}\di{}x=I_l & \text{on}\, E_l, \, 1\leq l\leq N_{\text{el}},\\
		\frac{1}{E_l}\int_{E_l}u\di{}x+z_l\int_{E_l}\bs{j}\cdot\bs{n}\di{}x=U_l & \text{on}\,E_l, \, 1\leq l \leq N_{\text{el}},\\
	\end{cases}
\end{equation}
where $\bs{n}$ is the unit outward normal. The first one is the no-flux condition, the second one means that the total injected current through each electrode is known and the third one means that the shared interface between electrode and material has an infinitesimally thin layer with surface impedance $z_l$. We have the following two constraints (Kirchhoff law of charge conservation and ground potential condition) to guarantee well-posedness:
\begin{equation}\label{eq:pde.bc.2}
	\sum_{l=1}^{N_{\text{el}}}I_l=0 \quad \text{and} \quad \sum_{l=1}^{N_{\text{el}}}U_l=0.
\end{equation}
The conductivity $\bs{\sigma}$ must be strictly positive to guarantee coercivity of the bilinear operator. The vector $\bs{I}=(I_1,\ldots,I_{N_{\text{el}}})$ provides deterministic injected currents, and the vector $\bs{U}=(U_1,\ldots,U_{N_{\text{el}}})$ provides random measurements of the potential at the electrodes. We solve the system \eqref{eq:pde.1},  \eqref{eq:pde.2} subject to the conditions \eqref{eq:pde.bc.1} and \eqref{eq:pde.bc.2} for the unknowns $(u, \bs{U})$, which are the potential field on $\cl{D}$ and the potential field at the electrodes, respectively.

\subsubsection{Experimental design formulation}

Ten electrodes are placed on the surface of the rectangular domain $\cl{D}=[0, 20]\times [0, 2]$, five on the top, and five on the bottom. The fiber orientations $\theta_1$ in ply $\cl{D}_1$ and $\theta_2$ in ply $\cl{D}_2$ are the parameters of interest with the following assumed prior distributions:
\begin{equation}
	\theta_1\sim\pi(\theta_1)=\cl{U}\left(-\frac{\pi}{4}-0.05, -\frac{\pi}{4}+0.05\right), \quad \theta_2\sim\pi(\theta_2)=\cl{U}\left(\frac{\pi}{4}-0.05, \frac{\pi}{4}+0.05\right).
\end{equation}
The orthotropic conductivity $\bs{\sigma}$ is given by
\begin{equation}\label{eq:nuis.dist}
	\sigma_j(\phi_i)=\exp(\mu_j+\phi_i), \quad j=1,2,3, \quad i=1,2,
\end{equation}
where $\phi_i=\sigma_{\phi}z_i$ is considered a nuisance parameter for $\bs{z}\sim\cl{N}(\bs{0},\bs{1})$,
corresponding to the two plies. We choose $\mu_1=\log(0.1)$ and $\mu_2=\mu_3=\log(0.02)$ in \eqref{eq:nuis.dist}.
Surface impedance is considered fixed at $z_l=z=0.1$ for each $1\leq l\leq N_{\rm{el}}$ to account for imperfections when attaching the electrodes so that a slight difference can exist between the potential fields at the electrodes. The data model is given by \eqref{eq:data.model} as follows:
\begin{equation}\label{eq:data.model.pde}
\bs{y}_i=\bs{g}_h(\bs{\theta},\bs{\phi})+\bs{\epsilon}_i, \quad \text{for}\, 1\leq i\leq N_e,
\end{equation}
where each $\bs{y}_i\in\bb{R}^q$, $q=N_{\text{el}}-1$ because one electrode is determined by \eqref{eq:pde.bc.2}. The error distribution is Gaussian, where $\bs{\epsilon}_i\sim\cl{N}(0, 10^{-4}\bs{1}_{q\times q}).$ We consider ten repetitions of the same experiment, so $N_e=10$. Moreover, $\bs{g}_h\equiv\bs{U}_h=(U_1,\ldots,U_{N_{\text{el}}-1})$ is given by a finite element approximation of $\bs{U}$ in the Galerkin sense. For details on the finite element formulation, see Appendix \ref{ap:fem} (see Figure \ref{fig:dom.grad}).

\begin{figure}[ht]
	\subfloat[Domain $\cl{D}$ consists of two plies\label{subfig-1:domain.ply}]{%
		\includegraphics[width=0.95\textwidth,trim={0 16cm 0 16cm},clip]{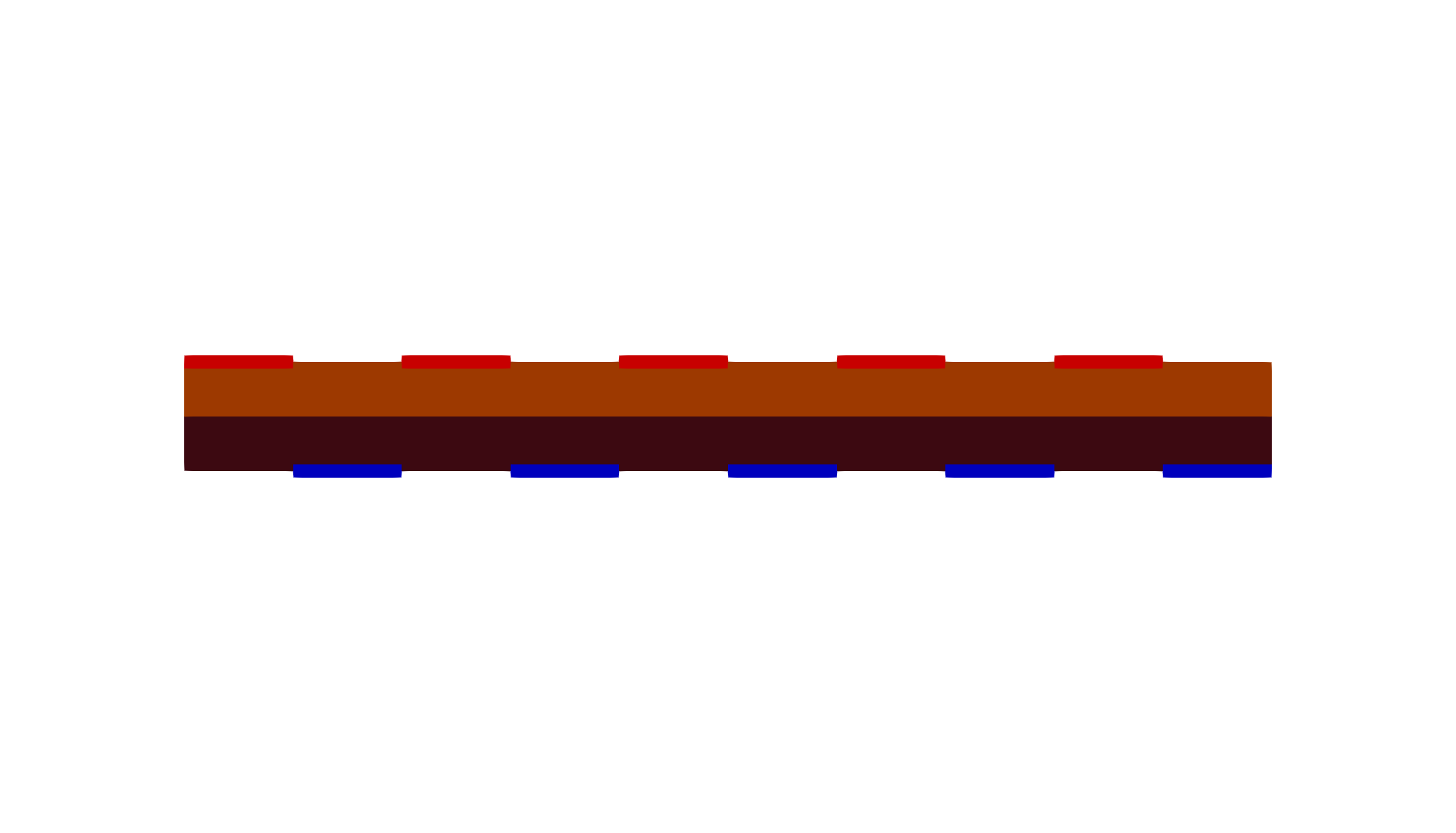}
	}
	\\
	\subfloat[Magnitude of the current and current streamlines\label{subfig-2:grad.ply}]{%
		\includegraphics[width=0.95\textwidth,trim={0 16cm 0 16cm},clip]{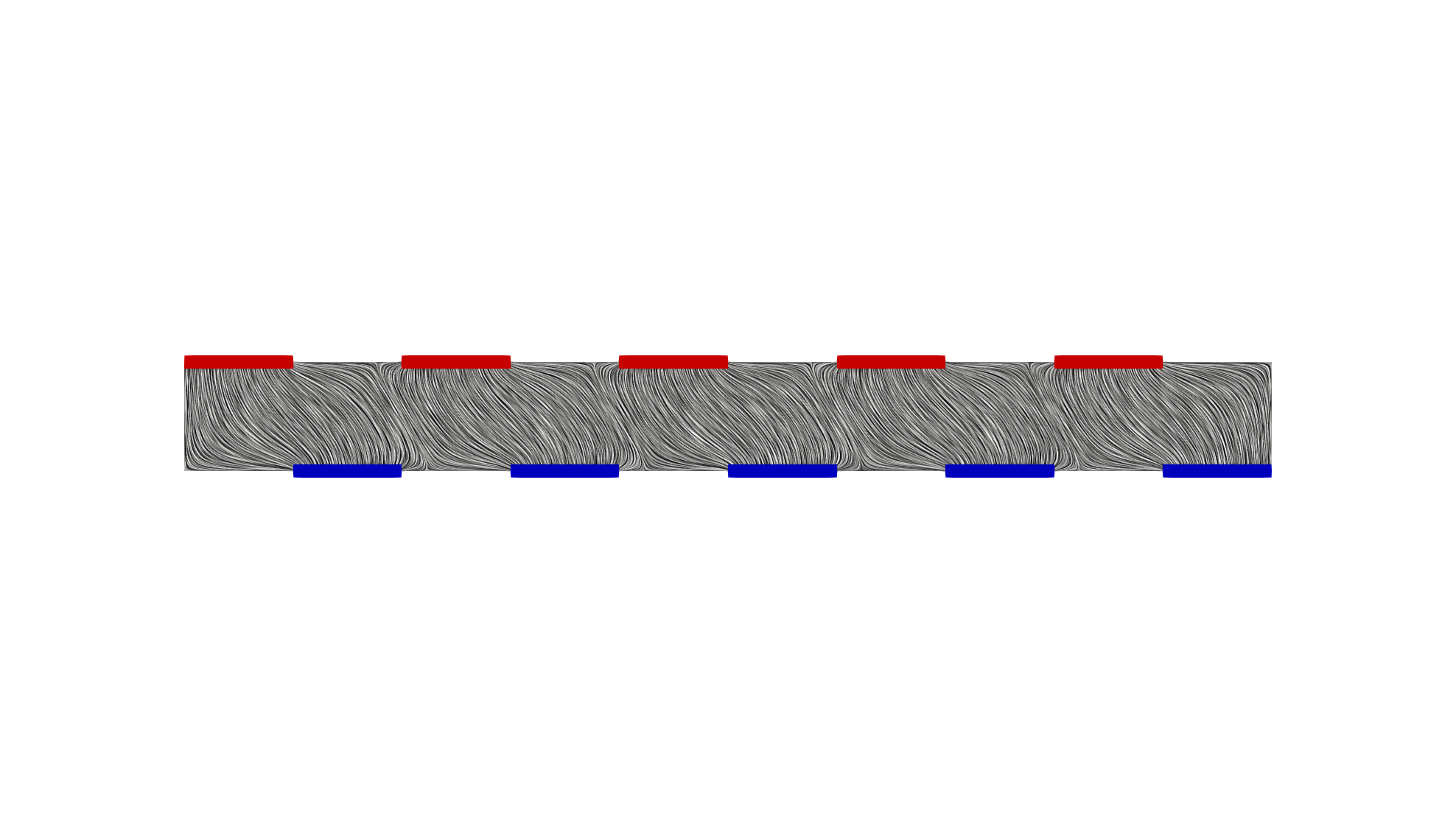}
	}
	\caption{Example 3: Panel \protect\subref{subfig-1:domain.ply} Domain $\cl{D}$. Panel \protect\subref{subfig-2:grad.ply} Current streamlines on top of the magnitude of the current. Inlet: red electrodes. Outlet: blue electrodes.}
	\label{fig:dom.grad}
\end{figure}

\subsubsection{Numerical results}

\begin{figure}[ht]
    \centering
    \includegraphics[width=0.6\textwidth]{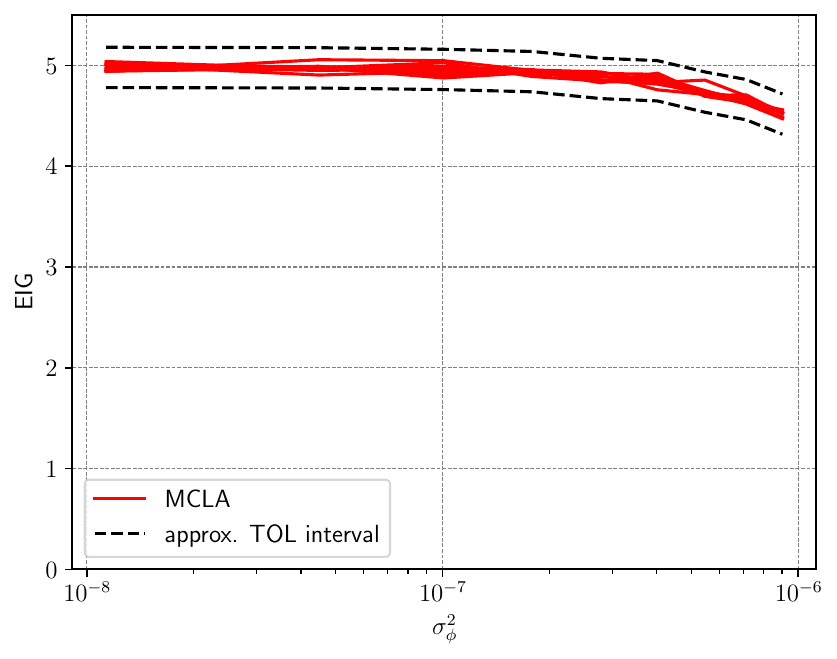}
    \caption{Example 3: EIG vs. $\sigma_{\phi}^2$ for ten runs of the MCLA estimator with approximate tolerance interval. Covariance matrix in \eqref{eq:small.likelihood.fam} loses positive definiteness when $\sigma_{\phi}^2$ exceeds a certain threshold of approximately $9.03\times10^{-7}$.}
    \label{fig:eit.bd}
\end{figure}

As mentioned, the small-noise approximation breaks down as the updated covariance matrix \eqref{eq:small.likelihood.fam} stops being positive definite for $\sigma_{\phi}^2$ larger than $9.03\times10^{-7}$ in the EIT example. Figure \ref{fig:eit.bd} illustrates ten runs of the MCLA estimator with the error tolerance specified at $TOL=0.2$. Evaluating the forward model $\bs{U}_h(\bs{\theta},\bs{\phi})$ for several values of $\bs{\theta}$ and $\bs{\phi}$ shows that results are approximately in the interval $[1,9]$ in absolute value. From \eqref{eq:data.model.pde}, we can see that $\bs{\epsilon}$ has standard deviation $10^{-2}$. The relative error in the data coming from the experimental noise is therefore between $0.1\%-1\%$. As for the norm of the gradient, we have $\lVert \nabla_{\bs{\phi}}\bs{g}(\bs{\theta},\bs{0})\rVert_2=\cl{O}(10)$, so Remark \ref{rmk:spc} suggests that the small-noise approximation is applicable for $\bs{\sigma}_{\bs{\phi}}\lessapprox 10^{-3}$, in accordance with the results presented in Figure \ref{fig:eit.bd}. Comparing the posterior covariance matrix for small and large nuisance uncertainty yields similar results as the ones presented in \ref{fig:eit.post.comp}, so we omitted them for brevity.

\begin{figure}[ht]
    \centering
    \includegraphics[width=0.6\textwidth]{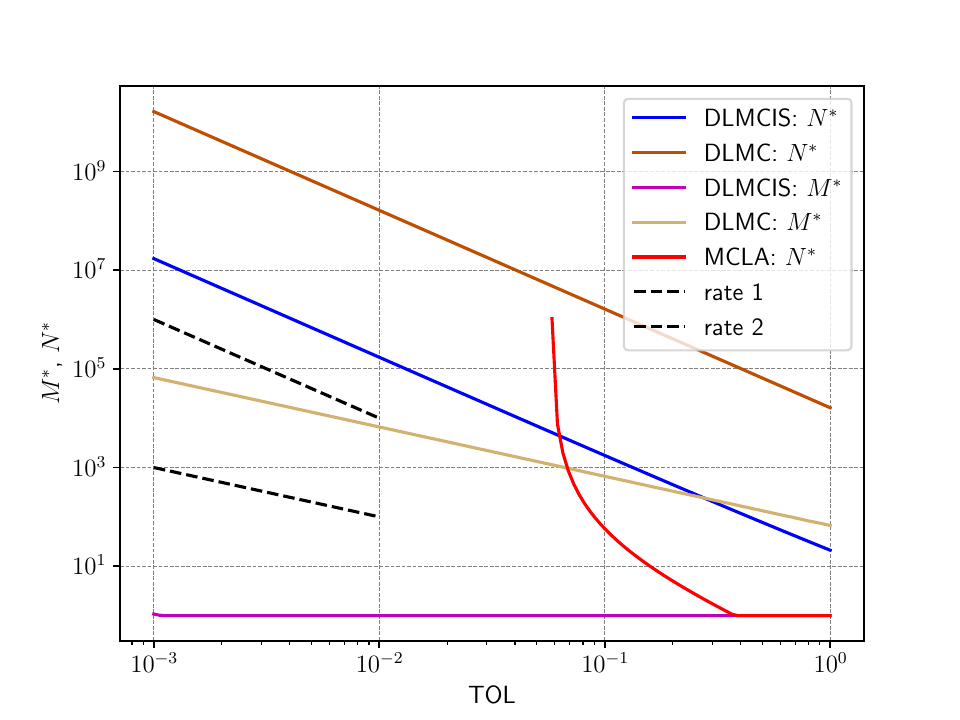}
    \caption{Example 3: Optimal number of outer ($N^{\ast}$) and inner ($M^{\ast}$) samples as a function of the tolerance $TOL$. The MCLA approximation cannot achieve accuracy when the inherent bias becomes greater than the desired tolerance, and $M^{\ast}$ for the DLMCIS estimator increases for small tolerances.}
    \label{fig:eit.opt}
\end{figure}

\begin{figure}[ht]
    \centering
    \includegraphics[width=0.6\textwidth]{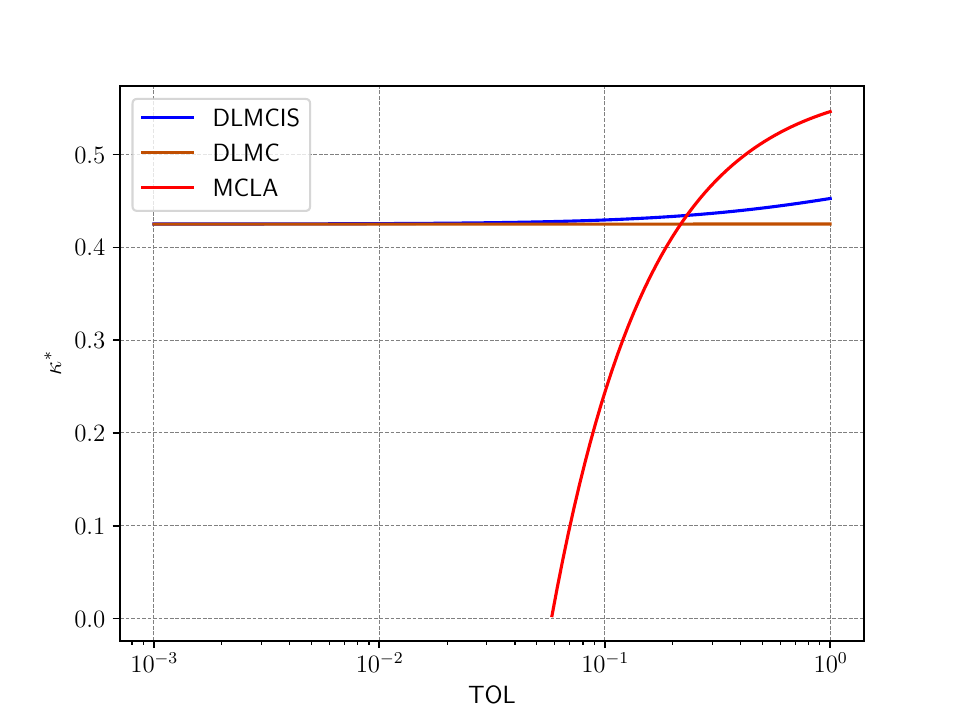}
    \caption{Example 3: The optimal splitting parameter $\kappa^{\ast}$ is constant for the DLMC and DLMCIS estimators for a small tolerance $TOL$. The MCLA estimator is infeasible for tolerances smaller than the inherent bias controlled by the number of repetitive experiments $N_e$, see \eqref{eq:laplace.bias}.}
    \label{fig:eit.nu}
\end{figure}

\begin{figure}[ht]
    \centering
    \includegraphics[width=0.6\textwidth]{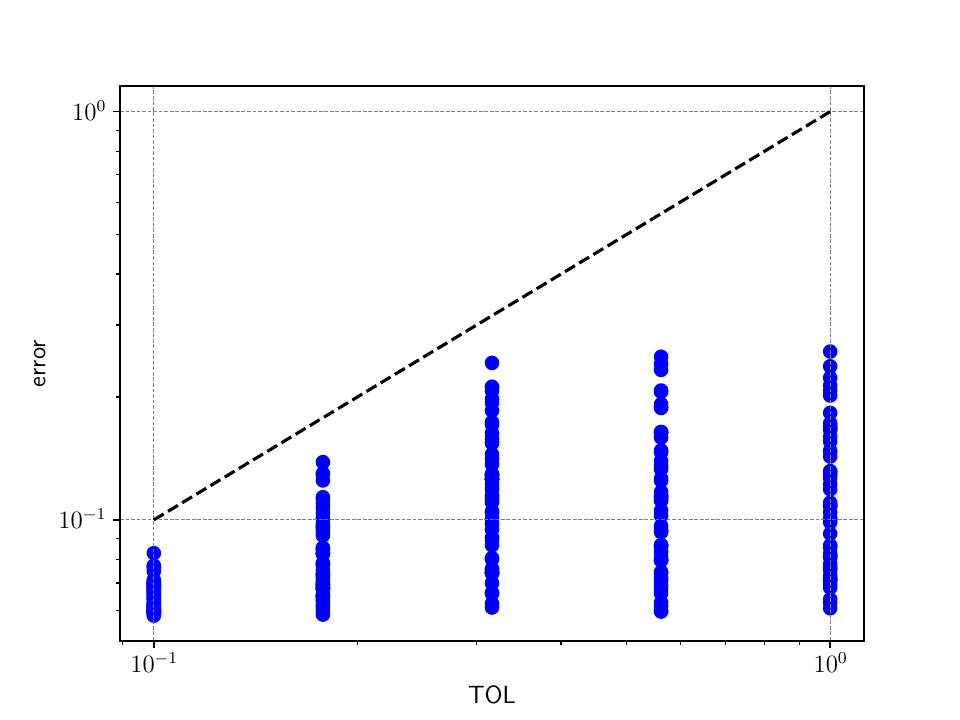}
    \caption{Example 3: Error vs. tolerance consistency plot. The MCLA estimator stays within the allowed tolerance $TOL$ with a predefined confidence parameter $C_\alpha=1.96$ ($\alpha=0.05$).}
    \label{fig:eit.evt}
\end{figure}

\begin{figure}[ht]
    \centering
    \includegraphics[width=0.6\textwidth,trim={0 1.5cm 0 1.5cm},clip]{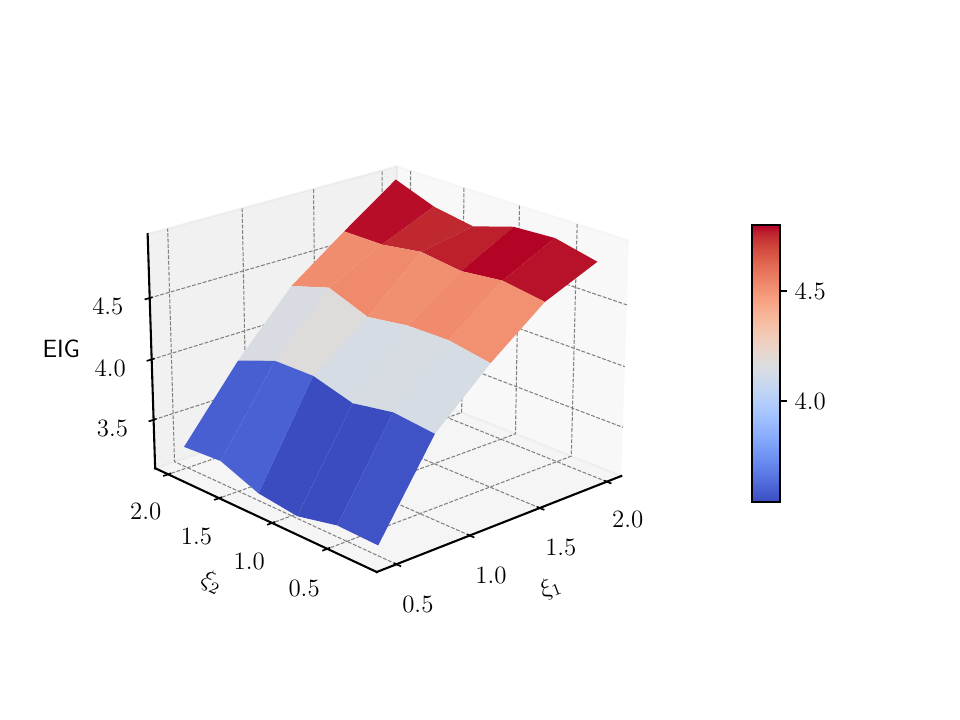}
    \caption{Example 3: Expected information gain as a function of the design parameters $\xi_1$ (the shift between top and bottom electrodes) and $\xi_2$ (the distance between electrodes) for the MCLA estimator.}
    \label{fig:eit.xi}
\end{figure}

Figure \ref{fig:eit.opt} reveals the optimal number of outer ($N^{\ast}$) and inner ($M^{\ast}$) samples for the DLMC, DLMCIS, and MCLA estimators for a fixed $\bs{\xi}=(2, 2)$, see Figure \ref{fig:eit.xi}, and $\sigma_{\phi}^2=10^{-10}$. The constants $C_{\mathrm{DL},1}$, $C_{\mathrm{DL},2}$ and $C_{\mathrm{DL},4}$ for the DLMC estimator were estimated using DLMC with $N=50$ and $M=30$. For the DLMCIS estimator, these constants were estimated using DLMCIS with $N=50$ and $M=10$ instead. Fewer inner samples are required in the latter case because of importance sampling. Finally,  the constants $C_{\mathrm{LA},1}$ and $C_{\mathrm{LA},2}$ for the MCLA estimator were estimated using MCLA with $N=50$. Because of importance sampling, the DLMCIS estimator only requires one inner sample in the optimal setting, except for small tolerances. The MCLA estimator cannot handle small tolerances as the inherent bias becomes greater than the total allowed tolerance, as in Figure \ref{fig:eit.nu}, where we demonstrate the splitting parameter between the bias and variance error. The error vs. tolerance consistency plot in Figure \ref{fig:eit.evt} indicates that the MCLA estimator is consistent and could even handle a coarser mesh. This likely results from an overestimation of the constants $\eta$ and $C_{\mathrm{LA},3}$ in \eqref{eq:mcla.h}.

In Figure \ref{fig:eit.xi}, we present the dependence of the EIG on the design parameters $\xi_1$ and $\xi_2$, signifying the shift between the top and bottom electrodes and the distance between the electrodes, respectively, via the MCLA estimator. The setup is optimal for both large $\xi_1$ and large $\xi_2$; however, the shift between electrodes seems to have a much greater influence on the EIG than the distance between electrodes.

\section{Conclusion}
We propose two estimators that incorporate nuisance uncertainty in the data model: the MCLA estimator based on a Laplace approximation of the inner integral and the DLMCIS estimator based on importance sampling. Both of these estimators use the small-noise approximation to deal with nuisance uncertainty in a computationally tractable manner. We found that these estimators can accurately approximate the EIG functional for analytical and physical examples as long as the nuisance uncertainty is small enough. This limitation is justified by only considering those parameters that are not primarily responsible for the overall uncertainty of the model as nuisance parameters.
Furthermore, we also demonstrated that including the small-noise approximation in the estimators does not affect their analytical error convergence derived in previous work \cite{Bec18}.

\section{Acknowledgments}

This publication is based upon work supported by the King Abdullah University of Science and Technology (KAUST) Office of Sponsored Research (OSR) under Award No.~OSR-2019-CRG8-4033, the Alexander von Humboldt Foundation, the Deutsche Forschungsgemeinschaft (DFG, German Research Foundation) -- 333849990/GRK2379 (IRTG Modern Inverse Problems), and was partially supported by the Flexible Interdisciplinary Research Collaboration Fund at the University of Nottingham Project ID 7466664.


\appendix

\section{Derivation of the work for two inner loops}\label{ap:double.inner.loop}
The average computational work of the DLMC estimator with two inner loops \eqref{EIG} is given by the solution of the optimization problem
\begin{equation}
  (N^{\ast},M_1^{\ast},M_2^{\ast})=\argmin_{(N,M_1,M_2)}N(M_1+M_2) \quad \text{subject to}\quad \frac{C_1}{M_1}+\frac{C_2}{M_2}\leq TOL,
\end{equation}
where $C_1=\bb{V}[p(\bs{Y}|\bs{\theta},\bs{\phi})|\bs{Y},\bs{\theta}]$, $C_2=\bb{V}[p(\bs{Y}|\bs{\theta},\bs{\phi})|\bs{Y}]$ and $TOL$ is the required tolerance. The derivation of the constants $C_1$ and $C_2$ follows from a simimlar argument as in \cite{Bec18}. Using the method of Lagrange multipliers, we obtain the equation
\begin{equation}
  \cl{L}(N,M_1,M_2)=N(M_1+M_2)-\lambda\left(\frac{C_1}{M_1}+\frac{C_2}{M_2}- TOL\right)
\end{equation}
with partial derivatives
\begin{equation}
  \frac{\partial\cl{L}}{\partial M_1}=N+\frac{\lambda C_1}{M_1^2}
\end{equation}
and
\begin{equation}
  \frac{\partial\cl{L}}{\partial M_2}=N+\frac{\lambda C_2}{M_2^2}.
\end{equation}
Setting these to zero immediately gives
\begin{equation}
M_2=\sqrt{\frac{C_2}{C_1}}M_1.
\end{equation}
From the partial derivative
\begin{equation}
  \frac{\partial\cl{L}}{\partial \lambda}=\frac{C_1}{M_1}+\frac{C_2}{M_2}- TOL,
\end{equation}
set to zero, we obtain that
\begin{equation}
  TOL=\frac{C_1+\sqrt{C_1C_2}}{M_1}.
\end{equation}
Thus, we have that
\begin{equation}
  M_1+M_2=\frac{C_1+2\sqrt{C_1C_2}+C_2}{TOL}.
\end{equation}

\section{Derivation of the order of the maximum a posteriori approximation}\label{ap:map}

Let
\begin{align}
  S(\bs{\theta})&=\frac{1}{2}N_e(\bs{g}(\bs{\theta}_t,\bs{0})-\bs{g}(\bs{\theta},\bs{0}))^{\trans}\bs{\Sigma}_{\bs{\epsilon}}^{-1/2}\bs{C}(\bs{\theta})\bs{\Sigma}_{\bs{\epsilon}}^{-1/2}(\bs{g}(\bs{\theta}_t,\bs{0})-\bs{g}(\bs{\theta},\bs{0}))\nonumber\\
  &+\sum_{i=1}^{N_e}\bs{\epsilon}_i\bs{\Sigma}_{\bs{\epsilon}}^{-1/2}\bs{C}(\bs{\theta})\bs{\Sigma}_{\bs{\epsilon}}^{-1/2}(\bs{g}(\bs{\theta}_t,\bs{0})-\bs{g}(\bs{\theta},\bs{0}))-h(\bs{\theta}) +D(\bs{\theta})
\end{align}

with gradient
\begin{align}
  \nabla_{\bs{\theta}}S(\bs{\theta})&=-N_e\nabla_{\bs{\theta}}\bs{g}(\bs{\theta},\bs{0})\bs{\Sigma}_{\bs{\epsilon}}^{-1/2}\bs{C}(\bs{\theta})\bs{\Sigma}_{\bs{\epsilon}}^{-1/2}(\bs{g}(\bs{\theta}_t,\bs{0})-\bs{g}(\bs{\theta},\bs{0}))\nonumber\\
  &+N_e(\bs{g}(\bs{\theta}_t,\bs{0})-\bs{g}(\bs{\theta},\bs{0}))^{\trans}\bs{\Sigma}_{\bs{\epsilon}}^{-1/2}\nabla_{\bs{\theta}}\bs{C}(\bs{\theta})\bs{\Sigma}_{\bs{\epsilon}}^{-1/2}(\bs{g}(\bs{\theta}_t,\bs{0})-\bs{g}(\bs{\theta},\bs{0}))\nonumber\\
  &-\sum_{i=1}^{N_e}\bs{\epsilon}_i\bs{\Sigma}_{\bs{\epsilon}}^{-1/2}\bs{C}(\bs{\theta})\bs{\Sigma}_{\bs{\epsilon}}^{-1/2}\nabla_{\bs{\theta}}\bs{g}(\bs{\theta},\bs{0})\nonumber\\
  &+\sum_{i=1}^{N_e}\bs{\epsilon}_i\bs{\Sigma}_{\bs{\epsilon}}^{-1/2}\nabla_{\bs{\theta}}\bs{C}(\bs{\theta})\bs{\Sigma}_{\bs{\epsilon}}^{-1/2}(\bs{g}(\bs{\theta}_t,\bs{0})-\bs{g}(\bs{\theta},\bs{0}))-\nabla_{\bs{\theta}}h(\bs{\theta}) +\nabla_{\bs{\theta}}D(\bs{\theta}),
\end{align}
where $\nabla_{\bs{\theta}}\bs{C}=-2\bs{\Sigma}_{\bs\epsilon}^{-1/2}\nabla_{\bs{\theta}}\nabla_{\bs{\phi}}\bs{g}(\bs{\theta},\bs{0})\bs{\Sigma}_{\bs{\phi}}^{1/2}\bs{M}^{\trans}\bs{\Sigma}_{\bs\epsilon}^{-1/2}$.

Applying Newton's method and $\nabla_{\bs{\theta}} S(\hat{\bs{\theta}})=\bs{0}$ implies that

$\nabla_{\bs\theta} S(\bs{\theta}_t)+\nabla_{\bs\theta}\nabla_{\bs\theta} S(\bs{\theta}_t)(\hat{\bs{\theta}}-\bs{\theta}_t)\approx\bs{0}$, so
\begin{equation}
	\bs{\hat{\theta}}\approx\bs{\theta}_t-\nabla_{\bs\theta}\nabla_{\bs\theta} S(\bs{\theta}_t)^{-1}\nabla_{\bs\theta} S(\bs{\theta}_t),
\end{equation}
providing
\begin{multline}
  \bs{\hat{\theta}}\approx\bs{\theta}_t-\Bigg(N_e\nabla_{\bs{\theta}}\bs{g}(\bs{\theta}_t,\bs{0})\bs{\Sigma}_{\bs{\epsilon}}^{-1/2}\bs{C}(\bs{\theta}_t)\bs{\Sigma}_{\bs{\epsilon}}^{-1/2}\nabla_{\bs{\theta}}\bs{g}(\bs{\theta}_t,\bs{0})-\sum_{i=1}^{N_e}\bs{\epsilon}_i\bs{\Sigma}_{\bs{\epsilon}}^{-1/2}\bs{C}(\bs{\theta}_t)\bs{\Sigma}_{\bs{\epsilon}}^{-1/2}\nabla_{\bs{\theta}}\nabla_{\bs{\theta}}\bs{g}(\bs{\theta}_t,\bs{0})\\
  -2\sum_{i=1}^{N_e}\bs{\epsilon}_i\bs{\Sigma}_{\bs{\epsilon}}^{-1/2}\nabla_{\bs{\theta}}\bs{C}(\bs{\theta}_t)\bs{\Sigma}_{\bs{\epsilon}}^{-1/2}\nabla_{\bs{\theta}}\bs{g}(\bs{\theta}_t,\bs{0})-\nabla_{\bs{\theta}}\nabla_{\bs{\theta}}h(\bs{\theta}_t)+\nabla_{\bs{\theta}}\nabla_{\bs{\theta}}D(\bs{\theta}_t)\Bigg)^{-1}\\
  \times\left(-\sum_{i=1}^{N_e}\bs{\epsilon}_i\bs{\Sigma}_{\bs{\epsilon}}^{-1/2}\bs{C}(\bs{\theta}_t)\bs{\Sigma}_{\bs{\epsilon}}^{-1/2}\nabla_{\bs{\theta}}\bs{g}(\bs{\theta}_t,\bs{0})-\nabla_{\bs{\theta}}h(\bs{\theta}_t)+\nabla_{\bs{\theta}}D(\bs{\theta}_t)\right).
\end{multline}
The leading order term in the denominator is of order $\cl{O}_{\bb{P}}\left(N_e\lVert\bs{\Sigma}_{\bs{\phi}}\rVert\right)$, and the leading order term in the numerator is of order $\cl{O}_{\bb{P}}\left(\sqrt{N_e}\lVert\bs{\Sigma}_{\bs{\phi}}\rVert\right)$, so we have $\bs{\hat{\theta}}=\bs{\theta}_t+\cl{O}_{\bb{P}}\left(\frac{1}{\sqrt{N_e}}\right)$ as $N_e\rightarrow\infty$.

\section{Derivation of the finite element formulation}\label{ap:fem}

Let $(\Omega,\cl{F},\bb{P})$ be a complete probability space with a set of outcomes $\Omega$, $\sigma$-field $\cl{F}$, and probability measure $\bb{P}$. We define $\cl{H}\coloneqq H^1(\cl{D}\times\bb{R}^{N_{\rm{el}}})$ as the space of the solution for the potential field $(u(\omega),\bs{U}(\omega))$ for $\omega\in\Omega$ and the Bochner space $L_{\bb{P}}^2(\Omega;\cl{H})$ as
\begin{equation}
  L_{\bb{P}}^2(\Omega;\cl{H})\coloneqq \left\{(u,\bs{U}):\Omega\to\cl{H} \quad \text{s.t.}\, \int_{\Omega}\lVert (u(\omega),\bs{U}(\omega))\rVert_{\cl{H}}^2\di{}\bb{P}(\omega)<\infty\right\}.
\end{equation}
For the bilinear form $B:\cl{H}\times\cl{H}\to\bb{R}$, given as
\begin{equation}
  B((u,\bs{U}),(v,\bs{V}))\coloneqq \int_{\cl{D}}\bs{j}\cdot \nabla v\di{}\cl{D}+\sum_{l=1}^{N_{\rm{el}}}\frac{1}{z_l}\int_{E_l}(\bs{U}_l-u)(\bs{V}_l-v)\di{}E_l,
\end{equation}
 we aim to determine $(u,\bs{U})\in L_{\bb{P}}^2(\Omega;\cl{H})$ such that the weak formulation
\begin{equation}
  \bb{E}[B((u,\bs{U}),(v,\bs{V}))]=\bs{I}_e\cdot\bb{E}[\bs{U}]
\end{equation}
is fulfilled for all $(v,\bs{V})\in L_{\bb{P}}^2(\Omega;\cl{H})$.


\footnotesize


\end{document}